%% file: main.tex
\documentclass{amsbook}		% Nicer than default article style:  less
\usepackage{amsmath,amsthm}     % Handy math stuff, theorem environments.
\usepackage{amssymb}            % Fancy math symbols.
\usepackage{euscript}           % Nice script font.
\usepackage[all,ps]{xy}            % XY-pic diagram package
\CompileMatrices		% speeds up XY-pic diagram compilation
\UseComputerModernTips		% use Knuth's computer modern fonts
\newcommand{\cat}[1]{\EuScript {#1}}

\newcommand{\lra}{\longrightarrow}

\newcommand{\ra}{\rightarrow}
\newcommand{\la}{\leftarrow}

\newcommand{\tdot}{^{^\centerdot}} 	% upper dot
\newcommand{\tdotdot}{^{^{\centerdot \centerdot}}} % upper double dot
\newcommand{\comp}{\circ}

\newcommand{\eps}{\epsilon}
\newcommand{\epso}{\overline \epsilon}
\newcommand{\nuo}{\overline \nu}

\newdir{ >}{{}*!/-5pt/\dir{>}}                  % Make better tailed arrows.
\newdir{ |}{{}*!/-5pt/\dir{|}}                  % Make better mapsto arrows.

\newtheorem{thm}{Theorem}[section]
\newtheorem{defn}[thm]{Definition}
\newtheorem{prop}[thm]{Proposition}

\newtheorem{cor}[thm]{Corollary}
\newtheorem{lemma}[thm]{Lemma}
\newtheorem{conj}[thm]{Conjecture}

\theoremstyle{definition}  % Bold headings and Roman body text.
\newtheorem{construction}[thm]{Construction}

\newtheorem{note}[thm]{Note}

\newcommand{\tens}              {\otimes}               %tensor
\DeclareMathOperator*{\colim}{colim\,}

\begin{document}

\mainmatter
\include{cover} 			
\include{abstract}
\include{intro}
\include{contents}
\include{ch1.1}
\include{ch1.2}
\include{ch1.3}
\include{chap2}

\appendix
\include{appa}
\include{appb}
\include{appc}

\include{biblio}
\end{document}

%% file: cover.tex
\bigskip
\bigskip
\bigskip

\begin{center}
{\huge {\bf Cofibrance and Completion}}

\bigskip

{\large by}

\bigskip

{\LARGE Andrei Radulescu-Banu}

\vspace{1.5in}

{\large Submitted to the Department of Mathematics}

{\large in partial fulfillment of the requirements}
{\large for the degree of}

{\large Doctor of Philosophy}

{\large at the}

{\large MASSACHUSETTS INSTITUTE OF TECHNOLOGY}

\bigskip

{\large February 1999}

\bigskip
\bigskip
\bigskip

\copyright 1998 Andrei Radulescu-Banu. All rights reserved.

\bigskip

The author hereby grants to MIT permission to reproduce

and to distribute publicly paper and electronic

copies of this thesis document in whole or in part.

\bigskip
\bigskip
\bigskip

Author: \hrulefill
\end{center}

\begin{flushright}
Department of Mathematics \\
November 30, 1998
\end{flushright}

\bigskip
\bigskip

\begin{center}
Certified by: \hrulefill
\end{center}

\begin{flushright}
Haynes R. Miller \\
Professor of Mathematics \\
Thesis Supervisor
\end{flushright}

\bigskip
\bigskip

\begin{center}
Accepted by: \hrulefill  
\end{center}

\begin{flushright}
Richard B. Melrose \\
Professor of Mathematics \\
Chairman, Departmental Commitee on Graduate Students
\end{flushright}

\begin{center}

\end{center}

%\title{Homotopy Triples and Cotriples, prel version}
%\author{by}
%\author{Andrei Radulescu-Banu}
%% Default to current date.
%\date{\today}
%\maketitle

%% file: abstract.tex
.
\newpage

\bigskip
\bigskip
\bigskip

\begin{center}
{\Large {\bf Cofibrance and Completion}}

\bigskip

by

\bigskip

{\large Andrei Radulescu-Banu}

\vspace{1in}

Submitted to the Department of Mathematics

on December 1, 1998 in Partial Fulfillment of the 

Requirements for the Degree of Doctor of Philosophy

\end{center}

\vspace{1in}

{\large { ABSTRACT}}

\bigskip

We construct the Bousfield-Kan completion with respect to a triple, for a model category.
In the pointed case, we construct a Bousfield-Kan spectral sequence that computes the relative homotopy groups of the completion of an object.

These constructions are based on the existence of a diagonal for the cofibrant-replacement functor constructed using the small object argument.

A central result that we use, due to Dwyer, Kan and Hirschhorn, is that in a model category homotopy limits commute with the function complex.

\vspace{1.5in}

Thesis Supervisor: Haynes R. Miller

Title: Professor of Mathematics

%% file: intro.tex
.
\newpage

\chapter*{Introduction}
\label{secIntro}
\bigskip
\bigskip

The Bousfield-Kan ${\mathbb Z}/p$ completion of simplicial sets is a very useful tool in Homotopy Theory.
Naturally enough, it is important to be able to mimick its construction in other contexts, for example for simplicial commutative algebras, simplicial groups, etc.

In these notes we generalize Bousfield and Kan's approach and construct, with respect to a triple, a completion functor in a model category.
In the case the model category is pointed, we also describe a Bousfield-Kan spectral sequence that computes the relative homotopy groups for the completion.

An important technical point in the constructions we perform is the existence of a diagonal for the cofibrant-replacement functor constructed using the small object argument.
In fact, we also show that if all acyclic fibrations are monomorphic then the small object argument can be changed (made ``less redundant'') so that the fibrant-replacement functor it constructs carries a codiagonal (i.e. a multiplication).

As for the general model category theory needed to understand these notes, we sytematically use that homotopy limits commute with the function complex (see Chap. \ref{chapHSS}, Sec. \ref{subsecHomLimFunComples}).
This result, due to Dwyer, Kan and Hirschhorn, allows us to prove general statements about model categories by just reducing them to statements about simplicial sets.

\bigskip
\bigskip

The classical approach to completion due to Bousfield and Kan will not immediately work for the case a model category.
To explain, let us recall the Bousfield-Kan completion of simplicial sets, and then, as an example, let us see what changes need to be made to perform the same construction for simplicial commutative algebras.

Fix some notations: $sSets$ denotes the category of simplicial sets, and $sSets_*$ the category of pointed simplicial sets.
Both of $sSets$, $sSets_*$ are model categories, in a standard way.
Also, let $R$ be a commutative ring, for example ${\mathbb Z}/p$ or ${\mathbb Z}$.

If $X. \in Ob(sSets_*)$ is a pointed simplicial set, denote by $RX.$ the free simplicial $R$-module on $X.$ , modulo the simplicial submodule generated by the basepoint.

The functor $R$ : $sSets_* \lra sR$-$mod$ is part of a triple ($R, \nu, M$)
\begin{center}
$\xymatrix{
	X. \ar[r]^\nu &
	RX. &
	R^2X. \ar[l]_M
}$
\end{center}
where $\nu$ is given by the unit of $R$, and $M$ by multiplication in $R$.

One forms the standard cosimplicial pointed space associated to the triple ($R, \nu, M$)

\begin{equation}
\label{eqnIntro1}
\xymatrix{ 
	RX. \ar[r] \ar@<1ex>[r] & 
	R^2X. \ar@<1ex>[l] \ar[r] \ar@<1ex>[r] \ar@<2ex>[r] & 
	R^3X. \;\; ... \ar@<1ex>[l] \ar@<2ex>[l]
	}
\end{equation}

\bigskip

The completion $X^{\wedge}_R$ of $X.$ is defined as $tot$ of the standard cosimplicial pointed space (\ref{eqnIntro1}).

Two things make this construction behave well from a homotopy theoretic point of view.
First, $RX.$ (and $R^n X.$) depend only on the pointed homotopy type of $X.$ 
Second, the standard cosimplicial pointed space (\ref{eqnIntro1}) is {\bf Reedy fibrant}, so $tot$ computes the homotopy limit of the cosimplicial pointed space, seen as a diagram of pointed spaces over the category $\Delta$.

\bigskip
\bigskip
\bigskip

Let us see what happens when we perform a similar construction for simplicial commutative algebras.
Details and a more general treatment of this example can be found in Chap. \ref{chap1}, Sec. \ref{secExamplesBKCompletion}.

We fix a field $k$ and denote by ${\cat A}_k$ the category of simplicial augmented $k$-algebras - ${\cat A}_k$ is a pointed model category, in a standard way.

Consider the adjoint pair ($Ab$, $i$) given by the pair abelianization in the category ${\cat A}_k$, embedding of the category $Ab{\cat A}_k$ of abelian objects of ${\cat A}_k$ into ${\cat A}_k$

\bigskip
\begin{center}
$Ab : {\cat A}_k \xymatrix{ \ar@<.5ex>[r] & \ar@<.5ex>[l] } Ab{\cat A}_k : i$
\end{center}
\bigskip

Denote $R$ = $i \comp Ab$ : ${\cat A}_k \lra {\cat A}_k$ . 
$R$ can be thought of as the abelianization functor, seen as a functor from the category ${\cat A}_k$ to itself.

$R$ is part of a triple ($R, \nu, M$), and for a simplicial algebra $X. \in Ob{\cat A}_k$ one forms as before the cosimplicial object (\ref{eqnIntro1}) (this time, it is a cosimplicial object in the category ${\cat A}_k$).

\bigskip

The crucial problem is, if $X. \ra Y.$ is a weak equivalence, then $RX. \ra RY.$ is not necessarily a weak equivalence as well.
Thus, the cosimplicial object (\ref{eqnIntro1}) does not have the right homotopy theoretic meaning anymore.
If $X.$, $Y.$ are {\bf cofibrant} and $X. \ra Y.$ is a weak equivalence, then $RX. \ra RY.$ is an equivalence.

\bigskip

We resolve this problem in the following way.

We take the cofibrant replacement functor $(S, \eps)$ constructed by the small object argument in ${\cat A}_k$

\begin{center}
$\xymatrix {
	SX \ar@{->>}[r]^{\eps}_{\sim} & X
}$
\end{center}
and give it a codiagonal, turning it into a cotriple ($S, \eps, \Delta$) (see Chap. \ref{chap1}, Sec. \ref{secCofRep}).

Then we mix this cotriple with the triple ($R, \nu, M$) to form the ``correct'' cosimplicial object

\begin{center}
\[ \xymatrix{ 
	SRSX \ar[r] \ar@<1ex>[r] & 
	SRSRSX \ar@<1ex>[l] \ar[r] \ar@<1ex>[r] \ar@<2ex>[r] & 
	SRSRSRSX \;\; ... \ar@<1ex>[l] \ar@<2ex>[l]
	}
\]
\end{center}
(see Chap. \ref{chap1}, Sec. \ref{secMixTriCot} for details on this construction).
We denote this cosimplicial object by ${\cat R}_S \tdot(X)$. 

\bigskip

We define the Bousfield-Kan completion of the algebra $X.$ just as an object in the homotopy category ${\bf ho}{\cat A}_k$, as the homotopy limit of the diagram ${\cat R}_S \tdot(X.)$

\begin{center}
$X^{\wedge}_{Ab} = {\mathbf R}\lim^{\cat \Delta} {\cat R}_S \tdot(X.)$
\end{center}

\bigskip
\bigskip
\bigskip

We have seen that the classical approach to construct the Bousfield-Kan completion needs to be adjusted for the case of a triple in a model category.
The same is true for the Bousfield-Kan spectral sequence.

In the case of completion of simplicial sets, the Bousfield-Kan spectral sequence is obtained as the ``classical'' homotopy spectral sequence of the cosimplicial pointed space given by \eqref{eqnIntro1}, and it computes the homotopy groups of the completion $\pi_* (X_R^\wedge)$.

For the example of simplicial algebras, we can construct a spectral sequence from scratch. In the cosimplicial object ${\cat R}_S \tdot(X.)$, forget the algebra structure to obtain a cosimplicial pointed space, and take the homotopy spectral sequence of that cosimplicial space.
This spectral sequence will compute the homotopy groups of the completion $\pi_* (X_{Ab}^\wedge)$.

\bigskip

Let us outline the construction of the Bousfield-Kan spectral sequence for the case of a triple in a pointed model category.

Let ${\cat M}$ be a pointed model category and $X \tdot$ a cosimplicial object in ${\cat M}$.
We define the homotopy spectral sequence of $X \tdot$ (see Chapter \ref{chapHSS}) as the ``classical'' homotopy spectral sequence of the cosimplicial pointed space ${\bf Hom}_*(W, X \tdot)$, where $W \in Ob{\cat M}$ is fixed.
${\bf Hom}_*($-$,$-$)$ denotes the {\it pointed} function complex (see Chap. \ref{chapHSS}, Sec. \ref{secPoiFunCom}).
The homotopy spectral sequence computes the relative homotopy groups with coefficients in $W$ of the homotopy limit ${\bf R}\lim^\Delta X \tdot$, that is,

\bigskip
\begin{center}
$\pi_* ({\bf Hom}_* (W, {\bf R}\lim^\Delta X \tdot))$
\end{center}
\bigskip

The Bousfield-Kan spectral sequence is defined as the homotopy spectral sequence of the cosimplicial object ${\cat R}_S \tdot(X)$.
Observe that the Bousfield-Kan spectral sequence computes the relative homotopy groups of the completion $X^{\wedge}$.

\bigskip
\bigskip
\bigskip
\bigskip
\bigskip

\subsection*{Organization of the paper}

Chapter \ref{chap1} describes, for a triple in model category, the construction of the Bousfield-Kan completion, and, in the case the model category is pointed, the construction of the associated Bousfield-Kan spectral sequence.

\bigskip

Chapter \ref{chapHSS} should be considered as a reference for general model category theoretic results needed throughout these notes, and it contains the construction of the homotopy spectral sequence of a cosimplicial object in a pointed model category.

\bigskip

Appendices \ref{app1}, \ref{app2} and \ref{app3} contain results needed in order to prove that the Bousfiel-Kan completion is independent of the choice of the cofibrant-replacement cotriple (as defined in Chap. \ref{chap1}, Sec. \ref{secSummUp}) used in its construction.

\bigskip
\bigskip
\bigskip

                % Local Variables: 
                % mode: latex 
                % tex-command: amslatex-command
                % tex-dvi-view-command: tex-oneside-view-command
                % tex-start-of-header: "\\special"
                % tex-end-of-header: "\\maketitle"
                % abbrev-mode: t
                % End:

%% file: contents.tex
\tableofcontents

%% file: ch1.1.tex
\chapter{The Bousfield-Kan Completion}
\label{chap1}

\bigskip
\bigskip
\bigskip

In this chapter we investigate the Bousfield-Kan completion with respect to a triple in a model category.
In the particular case when the model category is pointed we describe an associated Bousfield-Kan spectral sequence.

To be able to perform the Bousfield-Kan completion, the model category needs to have a cofibrant-replacement functor carrying a {\bf diagonal}.
Sections \ref{secCofRep} - \ref{secSummUp} investigate the existence of a diagonal for a cofibrant-replacement functor (and the existence of a codiagonal for a fibrant-replacement functor).

In detail, here is what Chapter \ref{chap1} contains:

\bigskip
\bigskip

In Section \ref{secCofRep}, we construct a diagonal for the cofibrant replacement functor ($S, \eps$)
\begin{center}
$\xymatrix {
	SX \ar@{->>}[r]^{\eps}_{\sim} & X
}$
\end{center}
constructed with the small object argument in a cofibrantly generated model category.
This diagonal turns ($S, \eps$) into a cotriple ($S, \eps, \Delta$). 

\bigskip
\bigskip

In Section \ref{secFibRep}, we deal with the fibrant replacement functor.
For a cofibrantly generated model category with the property that all acyclic cofibrations are mono- morphisms, we modify the small object argument to construct a fibrant replacement functor ($T, \nu$)
\begin{center}
$\xymatrix {
	X \ar@{>->}[r]^{\nu}_{\sim} & TX
}$
\end{center}
that admits a codiagonal (or multiplication), effectively turning ($T, \nu$) into a triple ($T, \nu, M$).
The proofs are not dual to those of Sec. \ref{secCofRep}, but the end result mirrors the one of Sec. \ref{secCofRep}.

\bigskip
\bigskip

In Section \ref{secSummUp}, we abstract the constructions in the previous two sections, and define cofibrant-replacement cotriples.

Bousfield-Kan completion can be constructed for model categories that have at least one cofibrant-replacement cotriple.
Cofibrantly generated model categories have cofibrant-replacement cotriples (by Sec. \ref{secCofRep}), and if all their acyclic cofibrations are monomorphic they also have fibrant-replacement triples (by Sec. \ref{secFibRep}).

\bigskip
\bigskip

Section \ref{secMixTriCot} describes a category theoretic construction, needed later for the construction of the Bousfield-Kan completion. 
For a category ${\cat C}$, given a cotriple ($S, \eps, \Delta$) and a triple ($R, \nu, M$),  we construct an augmented cosimplicial object in ${\cat C}$, by mixing the structure maps of our cotriple and triple
\begin{center}
\[ \xymatrix{ 
	{\overline {\cat R}}\; \tdot \!\!\! _S (X): &  
	SX \ar[r] & 
	SRSX \ar[r] \ar@<1ex>[r] & 
	SRSRSX \ar@<1ex>[l] \ar[r] \ar@<1ex>[r] \ar@<2ex>[r] & 
	SRSRSRSX \;\; ... \ar@<1ex>[l] \ar@<2ex>[l]
	}
\]
\end{center}
We denote ${\cat R}_S \tdot (X)$ the cosimplicial part of ${\overline {\cat R}}\; \tdot \!\!\! _S (X)$, obtained by dropping the augmentation $SX$.

\bigskip
\bigskip

In Section \ref{secBCCF}, we define the Bousfield-Kan completion with respect to a triple ($R, \nu, M$) in a model category ${\cat M}$.
For that, we assume that ${\cat M}$ has at least one cofibrant-replacement cotriple ($S, \eps, \Delta$), and we assume that $R$ takes weak equivalences between cofibrant objects in ${\cat M}$ to weak equivalences.

The Bousfield-Kan completion $X^{\wedge}_R$ of an object $X$ is defined as the homotopy limit of the cosimplicial diagram ${\cat R}_S \tdot(X)$ :

\bigskip
\begin{center}
$X^{\wedge}_R = {\mathbf R}\lim^{\cat \Delta} {\cat R}_S \tdot(X)$
\end{center}
\bigskip

We prove that the completion does not depend on the choice of the cofibrant-replacement cotriple ($S, \eps, \Delta$).

\bigskip
\bigskip

In Section \ref{subsecBKSpectrSeqCompl}, in the particular case when the model category is pointed, we construct the Bousfield-Kan spectral sequence, that computes the relative homotopy groups of the completion of an object.
The Bousfiled-Kan spectral sequence is obtained by just taking the homotopy spectral sequence (developed in Chapter \ref{chapHSS}) of the cosimplicial object ${\cat R}_S \tdot(X)$, and ultimately turns out to be from $E_2$ on independent on the choice of the cofibrant-replacement cotriple ($S, \eps, \Delta$).

\bigskip
\bigskip

In Section \ref{secExamplesBKCompletion}, as an example, we describe the completion of simplicial commutative algebras with respect to certain triples that arise form adjoint pairs of functors between simplicial algebras and simplicial modules.
This example will include the completion with respect to the abelianization of simplicial commutative algebras.

If $k$ is a field, we conjecture that connected simplicial augmented $k$-algebras are complete with respect to abelianization, and the associated {\bf absolute} Bousfield-Kan spectral sequence is convergent.

\bigskip
\bigskip
\bigskip

\section{A Diagonal for the Cofibrant-Replacement Functor}
\label{secCofRep}
\bigskip
\bigskip

In this Section, for the fibrant-replacement functor ($S, \eps$) constucted via the small object argument in a cofibrantly generated model category
\begin{center}
$\xymatrix {
	SX \ar@{->>}[r]^{\eps}_{\sim} & X
}$
\end{center}
we will construct a diagonal $\Delta$, which will turn the fibrant-replacement functor ($S, \eps$) into a cotriple ($S, \eps, \Delta$)
\begin{center} 
$ \xymatrix {
	S^2X &
	SX \ar[r]^\eps \ar[l]_\Delta &
	X
	} $
\end{center}

To construct $\Delta$, we need to understand the details of the small object argument.
Assume that ${\cat M}$ is a cofibrantly generated model category, and fix a generating set of cofibrations
\begin{center}
$\;\;\{ \xymatrix{ A_i \ar@{>->}[r]^{\phi_i} & B_i } \!\}_{i \eps I}$ 
\end{center}

Relative $I$-cell complexes are defined (Hirschhorn, \cite{Hirschhorn}) as transfinite compositions of pushouts of elements in 
\begin{center}
$\;\;\{ \xymatrix{ A_i \ar@{>->}[r]^{\phi_i} & B_i } \!\}_{i \eps I}$ 
\end{center}
(Hovey, \cite{Hovey} calls relative $I$-cell complexes ``regular $I$-cofibrations'')

Our assumptions about $\;\;\{ \xymatrix{ A_i \ar@{>->}[r]^{\phi_i} & B_i } \!\}_{i \eps I}$ are that:
\begin{list}{-}{}
 \item The class of retracts of relative $I$-cell complexes coincides with the class of cofibrations of ${\cat M}$. 
 \item There exists a regular ordinal $\lambda$ such that each $A_i$ is $\lambda$-small relative to the subcategory of cofibrations of ${\cat M}$.
\end{list}	

\bigskip

\subsection{The Small Object Argument}
Let us start by recalling the small object argument.
Denote $0$ the initial object in ${\cat M}$. 

The cofibrant replacement functor ($S, \eps$)
\begin{center}
$\xymatrix {
	SX \ar@{->>}[r]^{\eps}_{\sim} & X
}$
\end{center}
is constructed as 
\begin{center}
\[ \xymatrix {
0=S_0X \ar[r] \ar[dr]^{\eps_0} & S_1X \ar[r] \ar[d]^{\eps_1} & ... \ar[r] & S_{\beta}X \ar[r] \ar[dll]^{\eps_{\beta}} & ... \ar[r] & SX={\colim^{\beta < \lambda} S_{\beta}X} \ar[dllll]^{\eps = {\colim^{\beta < \lambda} {\eps_{\beta}}}}\\ 
& X & & & &
} \]
\end{center}

where:

\begin{list}{-}{}
 \item $\colim^{\beta < \lambda}$ are indexed by ordinals $\beta$, with $\beta < \lambda$
 \item If $\beta < \lambda$ is a limit ordinal, $S_{\beta}X={\colim^{\beta^{\prime} < \beta}} S_{\beta^{\prime}}X$ and ${\eps}_{\beta}={\colim^{\beta^{\prime} < \beta}} {\eps}_{\beta^{\prime}}$
 \item If $\beta + 1 < \lambda$ , then $\xymatrix{ S_{\beta}X \ar@{>->}[r] & S_{\beta + 1}X }$ is a pushout of a sum of copies of the generating cofibrations

\begin{center}
\[ \xymatrix {
	{\underset {{\cat D}_{\beta}(X)} \coprod} A_i \ar[r] \ar@{>->}[d]_{\coprod \phi_i} & S_{\beta}X \ar@{>->}[d] \ar[ddr]^{{\eps}_{\beta}} & \\
	{\underset {{\cat D}_{\beta}(X)} \coprod} B_i \ar[r] \ar[drr] & S_{\beta + 1}X \ar[dr]^{{\eps}_{\beta + 1}} & \\
	& & X
}
\]
\end{center}
\end{list}	

In the diagram above, ${\cat D}_{\beta}(X)$ is the set of all commutative diagrams of the form
\begin{center}
\[ \xymatrix {
	A_i \ar[r] \ar@{>->}[d]_{\phi_i} & S_{\beta}X \ar[d]^{{\eps}_{\beta}} \\
	B_i \ar[r] & X
}
\]
\end{center}
with $i \eps I$. 
We denote ${\tau}_{\beta , \beta^\prime} : S_\beta X \lra S_{\beta^\prime} X$ and ${\tau}_{\beta} : S_{\beta}X \lra SX$ the natural maps.

\bigskip

($S, \eps$) is indeed a cofibrant-replacement functor.
The object $SX$ is cofibrant, as a transfinite composition of cofibrations starting from the initial object. 
The natural map $\eps$ is an acyclic fibration, because it has the right lifting property with respect to the generating cofibrations
\begin{center}
$\;\;\{ \xymatrix{ A_i \ar@{>->}[r]^{\phi_i} & B_i } \!\}_{i \eps I}$ 
\end{center}

\bigskip

\subsection{Constructions needed for the diagonal}
Having constructed the cofibrant replacement functor ($S, \eps$), we want to construct on top of ($S, \eps$) a cotriple ($S, \eps, \Delta$).
In order to construct the diagonal $\Delta$, we have to go through a series of three intermediary constructions.

\begin{construction}
\label{cofCons1}
For any maps $a$, $b$ that make commutative the diagram with full arrows
\begin{center}
\[ \ \xymatrix@C=.5in{
	A_i \ar[r]^{a} \ar@{>->}[d]_{\phi_i}
	& S_{\beta}Y \ar[d]^{{\tau}_{\beta , \beta + 1}} \\
	B_i \ar@{-->}[r]^(.4){L_i^{\beta + 1}(a, b)} \ar[dr]_{b} 
	& S_{\beta + 1}Y \ar[d]^{\eps_{\beta + 1}} \\
	&
	Y
}
\]
\end{center}
construct the map $L_i^{\beta + 1}(a, b)$ as the composite

\begin{center}
$\xymatrix {
	B_i \ar[r]^-{inc_{a,b}} 
	& {\underset {{\cat D}_{\beta}(Y)} \coprod} B_j \ar[r]
	& S_{\beta + 1}Y 
}$
\end{center}
where $inc_{a,b}$ is the inclusion into the $a,b$-th factor.
\end{construction}

The map $L_i^{\beta + 1}(a, b)$ makes the extended diagram commutative, and it is natural in $Y$.

\bigskip

\begin{construction}
\label{cofCons2}
For any maps $a$, $b$ making commutative the diagram with full arrows
\begin{center}
\[ \ \xymatrix@C=.5in{
	S_{\beta}X \ar[r]^{a} \ar[d]_{{\tau}_{\beta , \beta + 1}} 
	& S_{\beta}Y \ar[d]^{{\tau}_{\beta , \beta + 1}} \\
	S_{\beta + 1}X \ar@{-->}[r]^{L_{\beta + 1}(a, b)} \ar[dr]_{b} 
	& S_{\beta + 1}Y \ar[d]^{\eps_{\beta + 1}} \\
	&
	Y
	}
\]
\end{center}
construct the map $L_{\beta + 1}(a, b)$ as follows.
$S_{\beta + 1}X$ is a pushout of $S_{\beta}X$, and many copies of $B_i$, $i \in I$.
Define $L_{\beta + 1}(a, b)$ on $S_{\beta}X$ as ${\tau}_{\beta , \beta + 1} a$, and define it on the copies of $B_i$ using Construction \ref{cofCons1}. 
\end{construction}
$L_{\beta + 1}(a, b)$ makes the extended diagram commutative, and it is natural in $X$ and $Y$.

\bigskip

Given any map $f : X \lra Y$, note that in the diagram
\begin{center}
\[ \ \xymatrix {
	S_{\beta}X \ar[r]^{S_{\beta}f} \ar[d]_{{\tau}_{\beta , \beta + 1}} 
	& S_{\beta}Y \ar[d]^{{\tau}_{\beta , \beta + 1}} \\
	S_{\beta + 1}X \ar@{-->}[r]^{S_{\beta + 1}f} \ar[dr]_{f \eps_{\beta + 1}}
	& S_{\beta + 1}Y \ar[d]^{\eps_{\beta + 1}} \\
	& Y
	}
\]
\end{center}
the following relation holds: $L_{\beta + 1}(S_{\beta}f, f \eps_{\beta + 1})$ = $S_{\beta + 1}f$.

Also, observe that in the diagram 
\begin{center}
\[ \ \xymatrix@C=.5in{
	S_{\beta}X \ar[r]^{a} \ar[d]_{{\tau}_{\beta , \beta + 1}} 
	& S_{\beta}Y \ar[r]^{a'} \ar[d]^{{\tau}_{\beta , \beta + 1}} 
	& S_{\beta}Z \ar[d]^{{\tau}_{\beta , \beta + 1}} \\
	S_{\beta + 1}X \ar@{-->}[r]^{L_{\beta + 1}(a, b)} \ar[dr]_{b} 
	& S_{\beta + 1}Y \ar@{-->}[r]^{L_{\beta + 1}(a', b')} \ar[d]^{\eps_{\beta + 1}} \ar[dr]_{b'}
	& S_{\beta + 1}Z \ar[d]^{\eps_{\beta + 1}} \\
	& Y
	& Z
	}
\]
\end{center}
we have the ``associativity'' relation
\begin{center}
$L_{\beta + 1}(a'a, b' L_{\beta + 1}(a, b))$ = $L_{\beta + 1}(a', b') L_{\beta + 1}(a, b)$
\end{center}

\bigskip

\begin{construction}
\label{cofCons4}
For any map $b$ : $SX \lra Y$, construct a lift $L(b)$
\begin{center}
\[ \ \xymatrix {
	SX \ar@{-->}[r]^{L(b)} \ar[dr]_{b} 
	& SY \ar[d]^\eps  \\
	& Y
	}
\]
\end{center}
as follows.
By transfinite induction on $\beta \leq \lambda$, for any map $b_\beta$ in the diagram
\begin{center}
\[ \ \xymatrix@C=.5in{
	S_\beta X \ar@{-->}[r]^{L_{\beta}^0 (b_\beta)} \ar[dr]_{b_\beta} 
	& S_\beta Y \ar[d]^{\eps_\beta} \\
	&
	Y
	}
\]
\end{center}
construct a map ${L_{\beta}^0 (b_\beta)}$, by taking $L_0^0(b_0)$ to be the identity map of the initial object, and using Construction \ref{cofCons2} at the succesive ordinal inductive step to construct ${L_{\beta + 1}^0 (b_{\beta + 1} )}$ = $L_{\beta + 1}(L_{\beta}^0 (b_{\beta + 1} \tau_{\beta, \beta + 1}), b_{\beta + 1})$.

Define $L(b)$ as $L_{\lambda}^0 (b)$, as a particular case of ${L_{\beta}^0}$ for $\beta$ = $\lambda$.
\end{construction}

\bigskip

The map $L(b)$ is natural in $X$ and $Y$, and satisfies 
\begin{equation}
\label{eqnunitalL0}
\eps L(b) = b
\end{equation}

Given any map $f : X \lra Y$, observe that in the diagram
\begin{center}
$ \ \xymatrix {
	SX \ar@{-->}[r]^{Sf} \ar[dr]_{f \eps}
	& SY \ar[d]^\eps \\
	& Y
	}
$
\end{center}
the following relation holds: 
\begin{equation}
\label{eqnunitalL}
L(f \eps) = Sf
\end{equation}
Also, in the diagram
\begin{center}
$ \ \xymatrix {
	SX \ar@{-->}[r]^{L(b)} \ar[dr]_{b} 
	& SY \ar[d]^\eps \ar@{-->}[r]^{L(b')} \ar[dr]_{b'}
	& SZ \ar[d]^\eps  \\
	& Y 
	& Z
	}
$
\end{center}
the associativity relation holds: 
\begin{equation}
\label{eqnassocL}
L(b' L(b)) = L(b') L(b)
\end{equation}

\bigskip

\subsection{The cotriple ($S, \eps, \Delta$)}
\label{subsecdelta}
We construct the diagonal $\Delta$ as $\Delta$ = $L(id_S)$, using Construction \ref{cofCons4} in the diagram 
\begin{center}
\[  \xymatrix {
	SX \ar@{-->}[r]^\Delta \ar[dr]_{id_S} 
	& S^2X \ar[d]^{\eps S} \\
	& SX
	}
\]
\end{center}

The cotriple axioms $\eps S \comp \Delta = id_S$, $S \eps \comp \Delta = id_S$ and $S \Delta \comp \Delta = \Delta S \comp \Delta$ are satisfied because of the two ``unital relations'' \eqref{eqnunitalL0}, \eqref{eqnunitalL} and the associativity relation \eqref{eqnassocL}.

\bigskip

It is worth to observe that giving ($S, \eps$) a cotriple structure ($S, \eps, \Delta$) is in fact {\it equivalent} to constructing $L($-$)$ satisfying the unital $\eps L(b) = b$, $L(f \eps) = Sf$ and associative $L(b' L(b)) = L(b') L(b)$ relations.
This is because given the cotriple ($S, \eps, \Delta$) we reconstruct $L($-$)$ by $L(b) = Sb \comp \Delta$, for $b : SX \lra Y$.

The equivalence between cotriple structures on ($S, \eps$) and unital and associative $L($-$)$ is related to the Kleisli category associated to a cotriple, cf. Mac Lane \cite{MacLane}, VI.5 Thm. 1.

\bigskip
\bigskip
\bigskip

\section{A Codiagonal for the Fibrant-Replacement Functor}
\label{secFibRep}
\bigskip
\bigskip

In the previous Section we showed that in a cofibrantly generated model category, the cofibrant-replacement functor constructed using the small object argument carries a cotriple structure.
It is natural to ask if the dual is true for the fibrant-replacement functor: does the fibrant-replacement functor constructed using the small object argument carry a {\it triple} structure?

The answer is no, but interestingly enough, if the model category is cofibrantly generated and {\it if all acyclic cofibrations are monomorphisms} then there is a variation (``less redundant'') of the small object argument that constructs a fibrant replacement functor that carries the structure of a triple.

\bigskip
\bigskip

Throughout this Section, ${\cat M}$ is a cofibrantly generated model category for which all acyclic cofibrations are monomorphisms.
We will start this Section with the description of the less redundant version of the small object argument, constructing a fibrant-replacement functor ($T, \nu$)
\begin{center}
$\xymatrix {
	X \ar@{>->}[r]^{\nu}_{\sim} & TX
}$
\end{center}
This construction is based on our extra hypothesis that all acyclic cofibrations are monomorphic.
In the second half of this Section we will construct on top of ($T, \nu$) a triple ($T, \nu, M$)
\begin{center}
$ \xymatrix {
	X \ar[r]^\nu &
	TX &
	T^2X \ar[l]_M
	} $
\end{center}

\bigskip
\bigskip

Fix a generating set of acyclic cofibrations 
\begin{center}
$\;\;\{ \xymatrix{ A_i \ar@{>->}[r]_{\sim}^{\phi_i} & B_i } \!\}_{i \eps I}$ 
\end{center}
(all are monomorphisms).

Our assumptions about $\;\;\{ \xymatrix{ A_i \ar@{>->}[r]_{\sim}^{\phi_i} & B_i } \!\}_{i \eps I}$ are now that:
\begin{list}{-}{}
 \item The class of retracts of relative $I$-cell complexes coincides with the class of acyclic cofibrations of ${\cat M}$. 
 \item There exists a regular ordinal $\lambda$ such that each $A_i$ is $\lambda$-small relative to the subcategory of acyclic cofibrations of ${\cat M}$.
\end{list}

\subsection{Small Object Argument, with less redundancy}
We construct ($T, \nu$) as the transfinite composition
\begin{center}
\[ \nu : \phantom{fo} \xymatrix@C=.3in{
X=T_0X \ar[r] & T_1X \ar[r] & ... \ar[r] & T_{\beta}X \ar[r] & ... \ar[r] & TX={\colim^{\beta < \lambda} T_{\beta}X}
}
\]
\end{center}

\bigskip

where:

\begin{list}{-}{}
 \item If $\beta < \lambda$ is a limit ordinal, $T_{\beta}X={\colim^{\beta^\prime < \beta}} T_{\beta^{\prime}}X$
 \item If $\beta + 1 < \lambda$ , then $\xymatrix {T_{\beta}X \ar@{>->}[r]_{\sim} & T_{\beta + 1}X}$ is a pushout of a sum of copies of the generating cofibrations

\begin{center}
\[ \xymatrix {
	{\underset {\cat D^\prime_\beta (X)} \coprod} A_i \ar[r] \ar@{>->}[d]_{\coprod {\phi_i}}^{\sim} & 
	T_{\beta}X \ar@{>->}[d]_{\sim} \\
	{\underset {\cat D^\prime_\beta (X)} \coprod} B_i \ar[r] & 
	T_{\beta + 1}X 
}
\]
\end{center}
The difference between our approach and the usual small object argument lies in what $\cat D^\prime_\beta (X)$ is.
In our case $\cat D^\prime_\beta (X)$ is the set of all maps $A_i \lra T_{\beta}X$ ($i \eps I$) that {\it do not admit a lift}
\begin{center}
$ \xymatrix {
	& T_{\beta^\prime}X \ar@{>->}[d]^\sim \\
	A_i \ar@{-->}[ur] \ar[r] & 
	T_\beta X
}$
\end{center}
for some $\beta^\prime < \beta$.
\end{list}	
 
\bigskip

Denote ${\nu}_{\beta}$ the maps $\xymatrix {{\nu}_{\beta} : X \ar@{>->}[r]_{\sim} & T_{\beta}X}$; we have $\nu = {\colim^{\beta < \lambda} {\nu}_{\beta}}$.

Denote $\tau_{\beta^\prime, \beta}$ the maps  $\tau_{\beta^\prime, \beta} : T_{\beta^\prime} X \lra T_\beta X$ , and $\tau_\beta$ the maps $\tau_\beta : T_\beta X \lra TX$.

\bigskip

To prove that ($T, \nu$) is a fibrant replacement functor, we first perform the following construction:

\begin{construction}
\label{consFib1}
For any map $a$ in the diagram  
\begin{center}
\[ \xymatrix@C=.5in{
	A_i \ar@{>->}[d]^\sim_{\phi_i} \ar[r]^a
	& T_\beta X \ar[d]^{\tau_{\beta, \beta + 1}}\\
	B_i \ar@{-->}[r]_(.4){L_i^{\beta + 1}(a)} 
	& T_{\beta + 1} X
	}
\]
\end{center}
we will construct a map $L_i^{\beta + 1}(a)$, so that the diagram is commutative.

If $a$ is in $\cat D^\prime_\beta (X)$, define $L_i^{\beta + 1}(a)$ as the composite
\begin{center}
$\xymatrix {
	B_i \ar[r]^-{inc_a} 
	& {\underset {{\cat D}_\beta^\prime(X)} \coprod} B_j \ar[r]
	& T_{\beta + 1}X 
}$
\end{center}
where $inc_a$ is the inclusion into the $a$-th factor.

If $a$ is not in $\cat D^\prime_\beta (X)$, it factors as
\begin{center}
\[ \xymatrix {
	& T_{\beta^\prime}X \ar@{>->}[d]^{\tau_{\beta^\prime, \beta}}_\sim \\
	A_i \ar@{-->}[ur]^{a'} \ar[r]_a & 
	T_\beta X
}
\]
\end{center}
for a minimal ordinal $\beta^\prime$, $\beta^\prime < \beta$.
Note that the map $a'$ is unique, because of our assumption that all acyclic cofibrations are monomorphic (so $\tau_{\beta^\prime, \beta}$ is a monomorphism).

We have that $a'$ $\eps$ $\cat D^\prime_{\beta^\prime} (X)$.
Define $L_i^{\beta + 1}(a)$ as the composite
\begin{center}
$\xymatrix@C=.5in{
	B_i \ar[r]^(.4){inc_{a'}} 
	& {\underset {{\cat D}^\prime_{\beta^\prime}(X)} \coprod} B_j \ar[r]
	& T_{\beta^\prime + 1}X \ar[r]^{\tau_{\beta^\prime + 1, \beta + 1}}
	& T_{\beta + 1}X
}$
\end{center}
where $inc_{a'}$ is the inclusion into the $a'$-th factor.
\end{construction}

\bigskip

We will give next a functor structure to $T_\beta$ (see Section \ref{subsecnat}).
Based on the functor structure of $T_\beta$, we will notice that the map $L_i^{\beta + 1}(a)$ of Construction \ref{consFib1} is natural in $X$ (see Lemma \ref{lemmaliftfib1}).
Finally, we will show that ($T$, $\nu$) is a fibrant-replacement functor.

\subsection{Functoriality}
\label{subsecnat}

Suppose we have a map
$\xymatrix{
X \ar[r]^f &
Y
}$.
For $\beta \leq \lambda$ we construct maps
$\xymatrix{
T_\beta X \ar[r]^{T_\beta f} 
& T_\beta Y
}$
by transfinite induction on $\beta$.

If $\beta$ is a limit ordinal, define $T_\beta f = {\colim^{\beta^\prime < \beta}} T_{\beta^{\prime}}f$.

If $T_\beta f$ was defined, we construct $T_{\beta + 1} f$ as follows:

\bigskip

For any map 
$\xymatrix{
A_i \ar[r]^-a 
& T_{\beta}X
}$ 
in $\cat D^\prime_\beta (X)$, Construction \ref{consFib1} in the diagram

\begin{center}
\[
\xymatrix{
	A_i \ar@{>->}[d]^\sim_{\phi_i} \ar[r]^-a &
	T_\beta X \ar[r]^{T_\beta f} &
	T_\beta Y \ar[d]^{\tau_{\beta, \beta + 1}} \\
	B_i \ar@{-->}[rr]_{L_i^{\beta + 1}(T_\beta f \comp a)} & &
	T_{\beta + 1} Y
}
\]
\end{center}
yields a map $L_i^{\beta + 1}(T_\beta f \comp a)$.
Recall that $T_{\beta + 1}X$ is defined as the pushout
\begin{center}
\[ \xymatrix {
	{\underset {\cat D^\prime_\beta (X)} \coprod} A_i \ar[r] \ar@{>->}[d]_{\coprod {\phi_i}}^{\sim} & 
	T_{\beta}X \ar@{>->}[d]_{\sim} \\
	{\underset {\cat D^\prime_\beta (X)} \coprod} B_i \ar[r] & 
	T_{\beta + 1}X 
}
\]
\end{center}
We define the map $T_{\beta + 1} f$ : $T_{\beta + 1} X \ra T_{\beta + 1} Y$ by glueing the two maps 

\bigskip
\begin{center}
$\tau_{\beta, \beta + 1} \comp T_\beta f$ : $T_\beta X \ra T_{\beta + 1} Y$

\bigskip

${\underset {\cat D^\prime_\beta (X)} \coprod} L_i^{\beta + 1}(T_\beta f \comp a)$ : ${\underset {\cat D^\prime_\beta (X)} \coprod} B_i \ra T_{\beta + 1}Y$
\end{center}
along ${\underset {\cat D^\prime_\beta (X)} \coprod} A_i$.
Observe that $T_{\beta + 1} f$ makes the diagram below commutative
\begin{center}
\[
\xymatrix{
	T_\beta X \ar[r]^{T_\beta f} \ar[d]_{\tau_{\beta, \beta + 1}} &
	T_\beta Y \ar[d]^{\tau_{\beta, \beta + 1}} \\
	T_{\beta + 1} X \ar[r]_{T_{\beta + 1} f} &
	T_{\beta + 1} Y
}
\]
\end{center}

\bigskip
\bigskip

This constructs $T_\beta$ on objects and on maps, for any ordinal $\beta$, $\beta \leq \lambda$. 
Not surprisingly, $T_\beta$ is a functor, as proved in the next Lemma:

\begin{lemma}
\label{lemmaFunctoriality}
In the construction above, $T_{\beta}$ : ${\cat M} \lra {\cat M}$ is a functor, and the maps $\nu_\beta$, $\nu$, $\tau_{\beta^\prime, \beta}$ and $\tau_\beta$ are natural, for any $\beta^\prime \leq \beta \leq \lambda$.

\end{lemma}

\begin{proof}
The essential step is to prove that $T_{\beta + 1}$ is a functor when we know
that $T_{\beta^\prime}$ are functors and that $\tau_{\beta^\prime, \beta}$ are natural maps, for all $\beta^\prime \leq \beta$.  

Given two maps
$\xymatrix{
X \ar[r]^f &
Y \ar[r]^g &
Z
}$
one has to show that $T_{\beta + 1}(gf)$ = $T_{\beta + 1}g \comp T_{\beta + 1}f$. 
This follows from how $T_{\beta^\prime}$, $\beta^\prime \leq \beta + 1$, are defined on maps.
We omit the details.
\end{proof}

\bigskip
\bigskip

The problem with the Construction \ref{consFib1} was that, at the time we defined the construction, we did not know that $T_{\beta}$ was a functor - because we used Construction \ref{consFib1} to constuct $T_{\beta}$ as a functor.
Since we now know the functor structure of $T_{\beta}$, we are ready to prove  naturality for the Construction  \ref{consFib1}:

\begin{lemma}
\label{lemmaliftfib1}
The map $L_i^{\beta + 1}(a)$ of Construction \ref{consFib1} is natural in $X$.
\end{lemma}

\begin{proof}
Consequence of how the functor $T_{\beta}$ is defined on maps, and of the fact that all acyclic cofibrations are monomorphic.
\end{proof}

\bigskip

At this stage, it is easy to show that the functor $T$ we constructed and the natural map $\nu$ form a fibrant replacement functor ($T$, $\nu$):

\begin{thm}
\label{thmCofRep}
The natural map
$\xymatrix {
	X \ar@{>->}[r]^-{\nu}_-{\sim} & TX
}$
is a trivial cofibration, and $TX$ is fibrant for any $X$.
\end{thm}

\begin{proof}
The map $\nu$ is an acyclic cofibration, as a transfinite composition of acyclic cofibrations.
$TX$ is fibrant for any $X$, because by Construction \ref{consFib2} below it has the right lifting property with respect to all the generating acyclic cofibrations 
\begin{center}
$\;\;\{ \xymatrix{ A_i \ar@{>->}[r]_{\sim}^{\phi_i} & B_i } \!\}_{i \eps I}$
\end{center}
\end{proof}

To complete the proof of Theorem \ref{thmCofRep} we perform
\begin{construction}
\label{consFib2}
For any map $a$ in the diagram 
\begin{center}
\[ \xymatrix {
	A_i \ar@{>->}[d]^\sim_{\phi_i} \ar[rd]^a
	& \\
	B_i \ar@{-->}[r]_-{L_i(a)}
	& TX
	}
\]
\end{center}
construct the map $L_i(a)$ as follows.

Since $A_i$ is $\lambda$-small relative to the subcategory of acyclic cofibrations of ${\cat M}$, the map $a$ factors as 
\begin{center}
\[\xymatrix {
	A_i \ar@{-->}[r]^-{a_\beta} \ar[dr]_a
	& T_\beta X \ar@{>->}[d]^{\tau_\beta}_\sim \\
	& TX
	}
\]
\end{center}
for a minimal ordinal $\beta$, $\beta < \lambda$. 
Observe that the map $a_\beta$ is unique, since $\tau_\beta$ is monomorphic.  

Using Construction \ref{consFib2}, we get a map $L_i^{\beta + 1}(a_\beta)$ : $B_i \ra T_{\beta + 1}X$.
We just define $L_i(a)$ as $\tau_{\beta + 1} L_i^{\beta + 1}(a_\beta)$ : $B_i \ra TX$ :
\begin{center}
\[ \xymatrix@C=.6in{
	A_i \ar@{>->}[d]^\sim_{\phi_i} \ar[r]^{a_\beta}
	& T_\beta (X) \ar@{>->}[d]_\sim^{\tau_{\beta, \beta + 1}} \\
	B_i \ar@{-->}[rd]_{L_i(a)} \ar[r]^-{L_i^{\beta + 1}(a_\beta)}
	& T_{\beta + 1} X \ar@{>->}[d]_\sim^{\tau_{\beta + 1}} \\
	& TX
} \]
\end{center}

\end{construction}

The map $L_i(a)$ just constructed is easily seen to be natural in $X$.

\bigskip

\subsection{Constructions needed for the codiagonal}
We build our way up to the construction of the triple ($T, \nu, M$), by doing two more constructions.

\begin{construction}
\label{consFib3}
For any map $a$ in the diagram
\begin{center}
\[ \xymatrix@C=.5in{
	T_\beta X \ar[d]_{\tau_{\beta, \beta + 1}} \ar[rd]^a
	& \\
	T_{\beta + 1} X \ar@{-->}[r]_-{L_{\beta + 1}(a)}
	& TY
	}
\]
\end{center}
we will construct a map $L_{\beta + 1}(a)$ that makes the diagram commutative.

$T_{\beta + 1} X$ is a pushout of $T_\beta X$ and of many copies of $B_i$, $i \in I$.
Define $L_{\beta + 1}(a)$ as $a$ on $T_\beta X$, and via Construction \ref{consFib2} for each copy of $B_i$ in the pushout.
\end{construction}

The map $L_{\beta + 1}(a)$ is natural in $X$, $Y$.
Given any map $f : X \lra Y$, in the diagram
\begin{center}
\[ \xymatrix {
	T_\beta X \ar[d]_{\tau_{\beta, \beta + 1}} \ar[r]^{T_\beta f}
	& T_\beta Y \ar[d]^{\tau_{\beta, \beta + 1}} \\
	T_{\beta + 1} X \ar@{-->}[dr]_{L_{\beta + 1}(a)} \ar[r]^{T_{\beta + 1} f}
	& T_{\beta + 1} Y \ar[d]^{\tau_{\beta + 1}} \\
	& TY
	}
\]
\end{center}
we have that $L_{\beta + 1}(\tau_\beta \comp T_\beta f)$ = $\tau_{\beta + 1} \comp T_{\beta + 1} f$.

\bigskip

\begin{construction}
\label{consFib5}
For any map $a$ in the diagram below, we will construct a map $L(a)$
\begin{center}
\[ \xymatrix {
	X \ar[d]_\nu \ar[dr]^a
	& \\
	TX \ar@{-->}[r]_{L(a)}
	& TY
	}
\]
\end{center}
that makes the diagram commutative.
To that end, for any $\beta \le \lambda$ and any map $a$ in the diagram 
\begin{center}
\[ \xymatrix {
	X \ar[d]_{\nu_\beta} \ar[dr]^a
	& \\
	T_\beta X \ar@{-->}[r]_{L_{\beta}^0 (a)}
	& TY
	}
\]
\end{center}
we construct a map ${L_{\beta}^0 (a)}$ that makes the diagram commutative, by transfinite induction on $\beta$.
We define $L_0^0(a)$ = $a$, and for the succesor ordinal induction step we use Construction \ref{consFib3} to define ${L_{\beta + 1}^0 (a)}$ = ${L_{\beta + 1} (L_\beta^0(a))}$.

We construct $L(a)$ = $L_{\lambda}^0 (a)$, as a particular case of the map $L_{\beta}^0$ for $\beta = \lambda$.
\end{construction}

The map $L(a)$ is natural in $X$ and $Y$, and satisfies
\begin{equation}
\label{eqnunitalLfib0}
L(a) \nu = a
\end{equation}

Given any map $f : X \lra Y$, in the diagram
\begin{center}
\[ \xymatrix {
	X \ar[d]_\nu \ar[rd]^{\nu f}
	& \\
	TX \ar@{-->}[r]_{Tf}
	& TY
	}
\]
\end{center}
it is easy to prove the relation

\begin{equation}
\label{eqnunitalLfib}
L(\nu f) = Tf
\end{equation}

\bigskip

The map $L(a)$ has another important property.
In the diagram 
\begin{center}
\[ \xymatrix {
	X \ar[d]_\nu \ar[dr]^a &
	Y \ar[d]_\nu \ar[dr]^{a'} & \\
	TX \ar@{-->}[r]_{L(a)} &
	TY \ar@{-->}[r]_{L({a'})} &
	TZ
	}
\]
\end{center}
we have that the associativity relation holds:

\begin{equation}
\label{eqnassocLfib}
L(L({a'}) a) = L({a'}) L(a)
\end{equation}

\bigskip

To prove associativity, one uses the following series of three intermediary results.
First, in the diagram 
\begin{center}
\[ \xymatrix{
	A_i \ar[rd]^a \ar@{>->}[d]^\sim_{\phi_i} &
	Y \ar[d]_\nu \ar[dr]^{a'} & \\
	B_i \ar@{-->}[r]_{L_i(a)} &
	TY \ar@{-->}[r]_{L({a'})} &
	TZ
	}
\]
\end{center}
we have that $L_i(L({a'}) a)$ = $L({a'}) L_i(a)$.
Second, in the diagram 
\begin{center}
\[ \xymatrix@C=.5in{
	T_\beta X \ar[rd]^a \ar[d]_{\tau_{\beta, \beta + 1}} &
	Y \ar[d]_\nu \ar[dr]^{a'} & \\
	T_{\beta + 1} X \ar@{-->}[r]_{L_{\beta + 1}(a)} &
	TY \ar@{-->}[r]_{L({a'})} &
	TZ
	}
\]
\end{center}
we have that $L_{\beta + 1}(L({a'}) a)$ = $L({a'}) L_{\beta + 1}(a)$.
Third, in the diagram 
\begin{center}
\[ \xymatrix{
	X \ar[rd]^a \ar[d]_{\nu_\beta} &
	Y \ar[d]_\nu \ar[dr]^{a'} & \\
	T_\beta X \ar@{-->}[r]_{L_{\beta}^0 (a)} &
	TY \ar@{-->}[r]_{L({a'})} &
	TZ
	}
\]
\end{center}
we have that $L_{\beta}^0(L({a'}) a)$ = $L({a'}) L_{\beta}^0(a)$.

\bigskip

\subsection{The triple ($T, \nu, M$)}

$M$ is constructed as $M$ = $L(id_T)$ in the diagram 
\begin{center}
\[ \xymatrix {
	TX \ar[d]_{\nu T} \ar[dr]^{id_T}
	& \\
	T^2X \ar@{-->}[r]_M
	& TX
	}
\]
\end{center}
The triple axioms for ($T, \nu, M$) are verified by arguments dual to those of Section \ref{subsecdelta}, using Construction \ref{consFib5} and the two unital \eqref{eqnunitalLfib}, \eqref{eqnunitalLfib0} and associativity \eqref{eqnassocLfib} relations.

\bigskip

Giving ($T, \nu$) a triple structure ($T, \nu, M$) is equivalent to constructing $L($-$)$ satisfying the unital $L(a) \nu = a$, $L(\nu f) = Tf$ and associative $L(L({a'}) a) = L({a'}) L(a)$ relations.
This is because from the triple ($T, \nu, M$) we can reconstruct $L($-$)$ by $L(a) = M \comp Ta$, for $a : X \lra TY$.

\bigskip
\bigskip
\bigskip
 
\section{Cofibrant-Replacement Cotriple. Fibrant-Replacement Triple.}
\label{secSummUp}
\bigskip
\bigskip

In order to construct the Bousfield-Kan completion with respect to a triple in a model category, we need to have at least a cotriple ($S, \eps, \Delta$) such that ($S, \eps$) is a cofibrant-replacement functor (we name such a cotriple ($S, \eps, \Delta$) a {\it cofibrant-replacement cotriple}).
This Section contains the necessary definitions, and sums up the results we proved in the previous two sections.

\bigskip
\bigskip

\begin{defn}
\label{defCofRepCot}
In a model category ${\cat M}$, a  cofibrant-replacement cotriple \\
($S, \eps, \Delta$) 
\begin{center} 
$ \xymatrix {
	S^2X &
	SX \ar[r]^{\eps}_{\sim} \ar[l]_\Delta &
	X
	} $
\end{center}
is a cotriple such that $\eps$ is a weak equivalence and $SX$ is cofibrant, for any $X$.
\end{defn}

\begin{list}{-}{Observe that:}
 \item We do not ask $\eps$ to be a fibration. 
 \item By the 2/3 axiom in ${\cat M}$, because of the cotriple structure, $\Delta$ is an equivalence as well.
 \item If all objects in ${\cat M}$ are cofibrant, the identity cotriple is a cofibrant-replace- \\
ment cotriple.
 \item In Section \ref{secCofRep}, if ${\cat M}$ is cofibrantly generated, we showed that the small object argument gives actually a cofibrant-replacement cotriple. 
A different choice of generating cofibrations yields generally a different cofibrant-replacement cotriple. 
 \item Dualizing the result of Section \ref{secFibRep}, if ${\cat M}$ is \emph{fibrantly} generated and its acyclic fibrations are epimorphisms, then cofibrant-replacement cotriples exist, constructed using the ``less redundant'' small object argument.
\end{list}

We summarize the results proved in Sections \ref{secCofRep}, \ref{secFibRep} as:
\begin{thm} 
\label{thmCofRepChar}
If ${\cat M}$ is either cofibrantly generated, or if all acyclic fibrations are epimorphisms and ${\cat M}$ is fibrantly generated, then cofibrant-replacement cotriples exist.
\end{thm}

The dual definition and theorem are:
\begin{defn}
\label{defFibRepTri}
In a model category ${\cat M}$, a  fibrant-replacement triple ($T, \nu, M$) 
\begin{center}
$ \xymatrix {
	X \ar[r]^{\nu}_{\sim} &
	TX &
	T^2X \ar[l]_M
	} $
\end{center}
is a triple such that $\nu$ is a weak equivalence and $TX$ is fibrant, for any $X$.
\end{defn}
\begin{thm} 
\label{thmFibRepChar}
If ${\cat M}$ is either fibrantly generated, or if all acyclic cofibrations are monomorphisms and ${\cat M}$ is cofibrantly generated, then fibrant-replacement triples exist.
\end{thm}

As an application of Section \ref{secFibRep}, note that $sSets$ has fibrant-replacement triples.

\bigskip
\bigskip
\bigskip

\section{The Mix of a Triple with a Cotriple}
\label{secMixTriCot}
\bigskip
\bigskip

This Section contains a triple-theoretic result needed for the construction of the Bousfield-Kan completion.
  
Let $\cat C$ be a category. Assume we have on $\cat C$ a cotriple ($S, \eps, \Delta$) and a triple ($R, \nu, M$):

\begin{center} 
$ \xymatrix {
	X &
	SX \ar[l]_\eps \ar[r]^\Delta &
	S^2X
	} $
	
$ \xymatrix {
	X \ar[r]^\nu &
	RX &
	R^2X \ar[l]_M
	} $
\end{center}

\bigskip
\bigskip

We will define below (Definition \ref{defcosres}) a cosimplicial resolution of $SX$
\begin{center}
$\xymatrix{ 
	{\overline {\cat R}}\; \tdot \!\!\! _S (X): &  
	SX \ar[r] &  
	{\cat R}_S \tdot (X)
}$
\end{center}
given by
\begin{center}
\[ \xymatrix{ 
	{\overline {\cat R}}\; \tdot \!\!\! _S (X): &  
	SX \ar[r] & 
	SRSX \ar[r] \ar@<1ex>[r] & 
	SRSRSX \ar@<1ex>[l] \ar[r] \ar@<1ex>[r] \ar@<2ex>[r] & 
	SRSRSRSX \;\; ... \ar@<1ex>[l] \ar@<2ex>[l]
	}
\]
\end{center}
natural in $X$, with the properties that 
\begin{equation}
\label{EqnContractionProperty1}
{\overline {\cat R}}\; \tdot \!\!\! _S (RX) {\rm \;\;\; contracts \; to \; the \; augmentation \;\;\;} SRX 
\end{equation}
and 
\begin{equation}
\label{EqnContractionProperty2}
R{\overline {\cat R}}\; \tdot \!\!\! _S (X) {\rm \;\;\; contracts \; to \; the \; augmentation \;\;\; } RSX
\end{equation}

\bigskip

The idea behind this construction is to exhibit the cosimplicial resolution ${\overline {\cat R}}\; \tdot \!\!\! _S (X)$ as the cosimplicial resolution of the triple associated to a certain adjoint pair.
The adjoint pair in question is going to be the (forget, free $S$-comodule) adjoint pair, composed with the (free $R$-module, forget) adjoint pair.

\bigskip
	
We will observe that if $f$, $g$ are any two triple maps from the triple ($R, \nu, M$) to another triple ($T, \nu, M$), then the maps of cosimplicial objects
\begin{center}
${\cat R}_S \,\!\!\!\! \tdot (X) \xymatrix{ \ar@<.5ex>[r]^f \ar@<-.5ex>[r]_g &} {\cat T}_S \,\!\!\! \tdot (X)$
\end{center}
are naturally homotopic (see Prop. \ref{propMapsOfMixedTriples}).

\bigskip
\bigskip

In this Section, we'll work with the whole resolution ${\overline {\cat R}}\; \tdot \!\!\! _S (X)$ rather than with the cosimplicial object ${\cat R}_S \tdot(X)$, because formulas are more compact for ${\overline {\cat R}}\; \tdot \!\!\! _S (X)$ .

We start with some preparations:

\bigskip
\bigskip

\subsection{Comodules over $S$ and modules over $R$.}  

Warning: in this Section we use slightly different-than-usual definitions for comodules (modules) over cotriples (resp. triples).

A comodule ($E, \Delta_E$) over the cotriple ($S, \eps, \Delta$) consists of a functor $E:{\cat C} \ra {\cat C}$ and a natural map
\begin{center} 
$ \xymatrix {
	E \ar[r]^{\Delta_E}&
	SE 
	} $
\end{center}
satisfying:
\bigskip

\begin{center} 
Coassociativity:
$ \xymatrix {
	E \ar[r]^{\Delta_E} \ar[d]_{\Delta_E} &
	SE \ar[d]^{\Delta S} \\
	SE \ar[r]^{S \Delta_E} &
	S^2E
	}$
\end{center}
\bigskip

\begin{center}
Counit:
$ \xymatrix {
	E &
	SE \ar[l]_{\eps E} \\
	&
	E \ar[ul]^{id} \ar[u]_{\Delta_E}
	}$ 
\end{center}

A comodule map is defined in the obvious way.

We denote $S$-$comod$ the category (in a higher universe) of comodules over the cotriple ($S, \eps, \Delta$).

\bigskip

Dually, a module ($F, M_F$) over the triple ($R, \nu, M$) satisfies by definition associativity and unit. 
We denote $R$-$mod$ the category of modules over ($R, \nu, M$).

\bigskip

\subsection{The extended comodule. The extended module.}

Denote $End{\cat C}$ the category of functors $F : {\cat C} \lra {\cat C}$, with natural maps as maps.

An example of a comodule ($E, \Delta_E$) over the cotriple ($S, \eps, \Delta$) is ($S, \Delta$) itself.

Even better, if $F : {\cat C} \lra {\cat C}$ is any functor, then ($SF, \Delta F)$ is a comodule over ($S, \eps, \Delta$) (the ``extended comodule'').
This association defines a right adjoint functor to the forgetful $S$-$comod \lra End{\cat C}$:

\bigskip

\begin{prop}
\label{propAdj1}
There is an adjoint pair ($\Phi_1=$forget$ , \Psi_1$) of functors
\begin{center}
$\Phi_1: \phantom{p} S$-$comod \xymatrix{ \ar@<.5ex>[r] & \ar@<.5ex>[l] } End{\cat C} \phantom{p} :\Psi_1$
\end{center}
where $\Phi_1 (E, \Delta_E) = E$ and $\Psi_1 (F) = (SF, \Delta F)$.
\end{prop}

Observe that $\Psi_1 (id_{\cat C})$ = ($S, \Delta$).
\begin{proof}
We construct adjointness morphisms. \\
First one: let comodule map $\phi : E \lra SF$ go to natural map
\begin{center}
$ \xymatrix { 
	\psi : E \ar[r]^\phi &
	SF \ar[r]^{\eps F} &
	F 
}$
\end{center}
Second one: let natural map $\psi : E \ra F$ go to comodule map 
\begin{center}
$ \xymatrix { 
	\phi : E \ar[r]^{\Delta_E} &
	SE \ar[r]^{S \psi} &
	SF 
}$ 
\end{center}
\end{proof}

\bigskip

Dually, we construct extended modules over the triple ($R, \nu, M$).
If $E : {\cat C} \lra {\cat C}$ is any functor, then the extended module is ($RE, ME$). 
This association defines a left adjoint functor to the forgetful functor $R$-$mod \lra End{\cat C}$:

\bigskip

\begin{prop}
\label{propAdj2}
There is an adjoint pair ($\Phi_2, \Psi_2=$forget) of functors
\begin{center}
$\Phi_2: \phantom{p} End{\cat C} \xymatrix{ \ar@<.5ex>[r] & \ar@<.5ex>[l] } R$-$mod \phantom{p} :\Psi_2$
\end{center}
where $\Phi_2 (E) = (RE, ME)$ and $\Psi_2 (F, M_F) = F$.
\end{prop}

\begin{proof}
Entirely dual to Proposition \ref{propAdj1}.
We construct adjointness morphisms. \\
First one: let natural map $\phi : E \lra F$ go to module map
\begin{center}
$ \xymatrix { 
	\psi : RE \ar[r]^{R \phi} &
	RF \ar[r]^{M_F} &
	F
}$ 
\end{center}
Second one:
let module map $\psi : RE \ra F$ go to natural map
\begin{center}
$ \xymatrix { 
	\phi : E \ar[r]^{\nu E} &
	RE \ar[r]^\psi &
	F 
}$ 
\end{center}
\end{proof}

\bigskip

Summing up Propositions \ref{propAdj1}, \ref{propAdj2}:
\begin{prop}
\label{propNiceAdjointness}

There is an adjoint pair ($\Phi, \Psi$) of functors
\begin{center}
$\Phi: \phantom{p} S$-$comod \xymatrix{ \ar@<.5ex>[r]^{\Phi_1} & \ar@<.5ex>[l]^{\Psi_1}} End {\cat C} \xymatrix{ \ar@<.5ex>[r]^{\Phi_2} & \ar@<.5ex>[l]^{\Psi_2} } R$-$mod \phantom{p} :\Psi$
\end{center}
\end{prop}

\bigskip
\bigskip

We have  $\Phi (E, \Delta_E) = (RE, ME)$, $\Psi (F, M_F) = (SF, \Delta F)$.

The composition of the adjunction morphisms in the proofs of Propositions \ref{propAdj1}, \ref{propAdj2} give adjunction morphisms for ($\Phi, \Psi$).
One is given by comodule map $\phi : E \lra SF$ going to module map 
\begin{center}
$ \xymatrix { 
	\psi : RE \ar[r]^{R\phi} &
	RSF \ar[r]^{R \eps F} &
	RF \ar[r]^{M_F} &
	F }$ 
\end{center}
and the other by module map $\psi : RE \lra F$ going to comodule map
\begin{center}
$ \xymatrix { 
	\phi : E \ar[r]^{\Delta_E} &
	SE \ar[r]^{S \nu E} &
	SRE \ar[r]^{S\psi} &
	SF }$ 
\end{center}

Denote these adjunction morphisms by $\nuo: id_{S-comod} \lra \Psi\Phi$, $\epso: \Phi\Psi \lra id_{R-mod}$.

\subsection{The cosimplicial resolution}
\label{subseccosres}

The adjoint pair ($\Phi, \Psi$) has as associated triple ($\Psi\Phi, \nuo, \Psi{\epso}\Phi$) on  $S$-$comod$.
This triple yields an augmented cosimplicial object ${\overline {\cat \Theta} \, \tdot}$ in the category $S$-$comod$
\begin{center}
\[
\xymatrix{ {\overline {\cat \Theta} \, \tdot} (E): &  E \ar[r] & \Psi\Phi(E) \ar[r] \ar@<1ex>[r] & \Psi\Phi\Psi\Phi(E) \ar@<1ex>[l] \ar[r] \ar@<1ex>[r] \ar@<2ex>[r] & \Psi\Phi\Psi\Phi\Psi\Phi(E) \;\; ... \ar@<1ex>[l] \ar@<2ex>[l]}
\]
\end{center}
\bigskip 

Denote $id_{\cat C}$ the identity functor of ${\cat C}$.

\bigskip

\begin{defn}
\label{defcosres}
Define ${\overline {\cat R}}\; \tdot \!\!\! _S = {\overline {\cat \Theta}} \, \tdot \Psi_1(id_{\cat C})$.
\end{defn}

Observe that ${\overline {\cat R}}\; \tdot \!\!\! _S R =  {\overline {\cat \Theta}} \, \tdot \Psi \Phi_2(id_{\cat C})$ and $R{\overline {\cat R}}\; \tdot \!\!\! _S = \Psi {\overline {\cat \Theta}} \, \tdot \Psi_1(id_{\cat C})$ are contractible.

\bigskip
\bigskip

In ${\overline {\cat R}}\; \tdot \!\!\! _S (X)$, coboundaries $d^i$ and codegeneracies $s^i$ are defined as follows:
\bigskip

\begin{center} 
${\overline {\cat R} \,}^n \!\!\! _S (X) = (SR)^{n+1}SX \ \ \ \ (n \ge -1)$

$d^i = (SR)^i [S \nu S \comp \Delta] (RS)^{n-i} : {\overline {\cat R} \,}^{n-1}\!\!\!\!\!\!\!\!\!\! _S \;\;\; (X) \lra {\overline {\cat R} \,}^n \!\!\! _S (X) \ \ \ \ (n \ge 0, 0 \le i \le n) $

$s^i = (SR)^i S [M \comp R \eps R ] S (RS)^{n-i} : {\overline {\cat R} \,}^{n+1} \!\!\!\!\!\!\!\!\!\! _S \;\;\; (X) \lra {\overline {\cat R} \,}^n \!\!\! _S (X) \ \ \ \ (n \ge 0, 0 \le i \le n) $
\end{center}
\bigskip

The contraction of ${\overline {\cat R}}\; \tdot \!\!\! _S (RY)$ to the augmentation $SRY$ is given by:

\begin{center}
$s^{n+1} = (SR)^{n+1} S [M \comp R \eps R ]: {\overline {\cat R} \,}^{n+1} \!\!\!\!\!\!\!\!\!\! _S \;\;\; (RY) \lra {\overline {\cat R} \,}^n \!\!\! _S (RY) \ \ \ \ (n \ge -1) $
\end{center}
and the contraction of $R{\overline {\cat R}} \; \tdot \!\!\! _S (X)$ to its augmentation $RSX$ is given by:
\begin{center}
$s^{-1} = [M \comp R \eps R ] S(RS)^{n+1} : R{\overline {\cat R} \,}^{n+1} \!\!\!\!\!\!\!\!\!\! _S \;\;\; (X) \lra R{\overline {\cat R} \,}^n \!\!\! _S (X) \ \ \ \ (n \ge -1) $
\end{center}
 
\bigskip

Using Prop. \ref{propNiceAdjointness} above and Appendix \ref{app1}, Prop. \ref{propMapsOfTriples} it is straightforward to prove the following Proposition:

\begin{prop}
\label{propMapsOfMixedTriples}
If $f$, $g$ are any two triple maps from a triple ($R, \nu, M$) to a triple ($T, \nu, M$), then the maps of cosimplicial objects
\begin{center}
${\cat R}_S \,\!\!\!\! \tdot (X) \xymatrix{ \ar@<.5ex>[r]^f \ar@<-.5ex>[r]_g &} {\cat T}_S \,\!\!\! \tdot (X)$
\end{center}
are naturally homotopic.
\end{prop}
For simplicity, we used the same notation $\nu$, $M$ for the structure maps of the triples ($R, \nu, M$) and ($T, \nu, M$).

\bigskip
\bigskip

\subsection{The simplicial resolution}
\label{subsecsimres}

We extend the results we obtained by duality. 
Suppose we have a triple ($T, \nu, M$) and a cotriple ($V, \eps, \Delta$) on $\cat C$:

\begin{center} 
$ \xymatrix {
	X \ar[r]^\nu &
	TX &
	T^2X \ar[l]_M
	} $

$ \xymatrix {
	X &
	VX \ar[l]_\eps \ar[r]^\Delta &
	V^2X
	} $
\end{center}

\bigskip

Then there exists a natural simplicial resolution of $TX$ 
\begin{center}
 \[ \xymatrix{
	{\overline {\cat V} \,}^T \!\!\! . (X): &
	{\cat V}^T \!\!\! .(X) \ar[r] &
	TX
	}
\]
\end{center}
given by
\begin{center}
 \[ \xymatrix{
	{\overline {\cat V} \,}^T \!\!\! . (X): & 
	... \;\; TVTVTVTX\ar[r] \ar@<1ex>[r] \ar@<2ex>[r] &
	TVTVTX \ar[r] \ar@<1ex>[r] \ar@<1ex>[l] \ar@<2ex>[l] &
	TVTX \ar[r] \ar@<1ex>[l] &
	TX
	}
\]
\end{center}
with ${\overline {\cat V} \,}^T \!\!\! . (TX)$, $T({\overline {\cat V} \,}^T \!\!\! . (X))$ contractible.
The formulas for $d_i$'s and $s_i$'s for ${\overline {\cat V} \,}^T \!\!\! . (X)$, ${\overline {\cat V} \,}^T \!\!\! . (TX)$ and $T({\overline {\cat V} \,}^T \!\!\! . (X))$ are dual to those for $d^i$, $s^i$.

\bigskip

Dual to Proposition \ref{propMapsOfMixedTriples} we have

\begin{prop}
\label{propMapsOfMixedCotriples}
If $f$, $g$ are two cotriple maps from a cotriple ($U, \eps, \Delta$) to a cotriple ($V, \eps, \Delta$), then the maps of simplicial objects
\begin{center}
${\overline {\cat U} \,}^T \!\!\! . (X) \xymatrix{ \ar@<.5ex>[r]^f \ar@<-.5ex>[r]_g &} {\overline {\cat V} \,}^T \!\!\! . (X)$
\end{center}
are naturally homotopic.
\end{prop}

\bigskip
\bigskip
\bigskip

                % Local Variables: 
                % mode: latex 
                % tex-command: amslatex-command
                % tex-dvi-view-command: tex-oneside-view-command
                % tex-start-of-header: "\\special"
                % tex-end-of-header: "\\maketitle"
                % abbrev-mode: t
                % End:

%% file: ch1.2.tex
\section{The Bousfield-Kan Completion Functor}
\label{secBCCF}
\bigskip
\bigskip
\bigskip

In this Section, for a model category ${\cat M}$ and a triple ($R, \nu, M$) (satisfying certain conditions, cf. Section \ref{subsecConstOfCompl}) we construct a Bousfield-Kan completion functor 
\begin{center}
\[({\rm -})^{\wedge}_R : {\bf ho} {\cat M} \lra {\bf ho} {\cat M}\]
\[X \mapsto X^{\wedge}_R\]
\end{center}
and we prove some of its properties.
For the particular case of a {\bf pointed} model category, we develop a Bousfield-Kan spectral sequence that computes the relative homotopy groups of the completion of an object.

We prove that the Bousfield-Kan completion is independent of the choice of cofibrant-replacement cotriple ($S, \eps, \Delta$) used in their definition.
This essentially relies on a repeated use of Appendix \ref{app1}, Thm. \ref{thmHomot2}, which says that in a model category two homotopic maps between cosimplicial objects 

\begin{center}
$f$ $\simeq$ $g$ : $X \tdot \lra Y \tdot$ 
\end{center}
induce equal maps on homotopy limits: ${\mathbf R}\lim^\Delta f$ = ${\mathbf R}\lim^\Delta g$.

We show that the Bousfield-Kan spectral sequence is independent of the choice of cofibrant-replacement cotriple from $E_2$ on.

\bigskip

We prove that completion is natural in the triple ($R, \nu, M$) in the following strong sense.
For two maps of triples $f$, $g$ from a triple ($R, \nu, M$) to a triple ($T, \nu, M$) in ${\cat M}$, we have

\begin{center}
$f^\wedge$ = $g^\wedge$ : $X^{\wedge}_R \lra X^{\wedge}_T$
\end{center}

\bigskip

The material in this section depends on the previous sections, but also on general theory of ${\bf R}\lim$ developed in Chapter \ref{chapHSS}, as well as on the technical results of Appendices \ref{app1}, \ref{app2} and \ref{app3}.

One could of course dualize all the statements in this Section, but we won't formulate the dual results for brevity.

\bigskip

\subsection{Construction}
\label{subsecConstOfCompl}
Let ${\cat M}$ be a model category. 
We fix a cofibrant replacement cotriple ($S, \eps, \Delta$) in ${\cat M}$.
We assume ${\cat M}$ has one; this is true for example if ${\cat M}$ satisfies one of the hypotheses of Theorem \ref{thmCofRepChar}.

\bigskip

Suppose ($R, \nu, M$) is a triple on ${\cat M}$.
We assume that for all equivalences $f$ between cofibrant objects of ${\cat M}$, $R(f)$ is an equivalence.

Using the cotriple ($S, \eps, \Delta$) and the triple ($R, \nu, M$), we form as in Section \ref{secMixTriCot} the cosimplicial resolution
\begin{center}
$ \xymatrix{ 
	{\overline {\cat R}}\; \tdot \!\!\! _S (X): &  
	SX \ar[r] &  
	{\cat R}_S \tdot (X)
}$
\end{center}

\begin{defn}
For $X$ an object in ${\cat M}$, the Bousfield-Kan $R$-completion $X^{\wedge}_R$ is defined as
\begin{center}
$X^{\wedge}_R = {\mathbf R}\lim^{\cat \Delta} {\cat R}_S \tdot(X)$
\end{center}
\end{defn}

\bigskip
\bigskip

If no danger of confusion, we will refer to the $R$-completion $X^{\wedge}_R$ as $X^{\wedge}$.

It is easy to observe that $X^{\wedge}$ depends only on the weak equivalence class of $X$.

The augmentation of ${\overline {\cat R}}\; \tdot \!\!\! _S (X)$ yields in ${\bf ho}{\cat M}$ a completion map 

\bigskip

\begin{center}
$X \lra X^{\wedge}$
\end{center}

\begin{thm}
\label{thmBCdef}
Completion $(-)^{\wedge}$ : ${\bf ho}{\cat M} \ra {\bf ho}{\cat M}$ is a functor.
The completion map $X \lra X^{\wedge}$ is a natural map.

The completion functor $(-)^{\wedge}$ and the completion map do not depend on the choice of cofibrant-replacement cotriple ($S, \eps, \Delta$) used in their definition.
\end{thm}

\begin{proof}
The difficult part is proving the second paragraph, so we will just describe that part of the proof.

Let ($S', \eps, \Delta$) be a second cofibrant-replacement cotriple - we denote its structure maps again $\eps$, $\Delta$, as for ($S, \eps, \Delta$).

Denote $X \lra X^{\wedge}_{R,S}$ the completion map constructed using the cotriple ($S, \eps, \Delta$), and $X \lra X^{\wedge}_{R,S'}$ the one constructed using ($S', \eps, \Delta$).

Using a double complex argument, we will construct in ${\bf ho}{\cat M}$ a natural isomorphism 
\begin{center}
$ \xymatrix{
	X \ar[r] \ar[dr] &
	X^{\wedge}_{R,S} \ar[d]^{\theta_{S', S}}_\cong \\
	& X^{\wedge}_{R,S'}
} $
\end{center}
with the property that $\theta_{S'', S'}\theta_{S', S}$ = $\theta_{S'', S}$ and $\theta_{S, S}$ = $id$.

To construct $\theta_{S', S}$, consider the following commutative diagram:

\bigskip

\begin{equation}
\label{eqnDoubleComplex}
\xymatrix{ 
	S'SX \ar[r] \ar[d] &
	S'SRSX \ar[r] \ar@<1ex>[r] \ar[d] & 
	S'S(RS)^2X \ar@<1ex>[l] \ar[d] & 
	... \\ 
	S'RS'SX \ar[r] \ar[d] \ar@<1ex>[d] & 
	S'RS'SRSX \ar[r] \ar@<1ex>[r] \ar[d] \ar@<1ex>[d] & 
	S'RS'S(RS)^2X \ar@<1ex>[l] \ar[d] \ar@<1ex>[d] & 
	... \\
	(S'R)^2S'SX \ar[r] \ar[d] \ar@<1ex>[d] \ar@<2ex>[d] \ar@<1ex>[u] & 
	(S'R)^2S'SRSX \ar[r] \ar@<1ex>[r] \ar[d] \ar@<1ex>[d] \ar@<2ex>[d] \ar@<1ex>[u] & 
	(S'R)^2S'S(RS)^2X \ar@<1ex>[l] \ar[d] \ar@<1ex>[d] \ar@<2ex>[d] \ar@<1ex>[u] &
	... \\
	(S'R)^3S'SX \ar[r] \ar@<1ex>[u] \ar@<2ex>[u] & 
	(S'R)^3S'SRSX \ar[r] \ar@<1ex>[r] \ar@<1ex>[u] \ar@<2ex>[u] & 
	(S'R)^3S'S(RS)^2X \ar@<1ex>[l] \ar@<1ex>[u] \ar@<2ex>[u] & 
	... \\
	... & ... & ... & ...
	}
\end{equation}

\bigskip

Observe that the top row in \eqref{eqnDoubleComplex} is $S'{\overline {\cat R}}\; \tdot \!\!\! _{S} (X)$, and the left column is ${\overline {\cat R}}\; \tdot \!\!\! _{S'} (SX)$.

\bigskip

Let us denote by ${\overline {\cat R}}\; \tdotdot \!\!\!\!\! _{S',S} (X)$ the diagram \eqref{eqnDoubleComplex}.
Denote by ${\cat R} \tdotdot \!\!\! _{S',S} (X)$ the subdiagram obtained by erasing in \eqref{eqnDoubleComplex} the top row and the left column.

\bigskip

${\cat R} \tdotdot \!\!\! _{S',S} (X)$ is a bicosimplicial diagram, and the diagram \eqref{eqnDoubleComplex} yields cosimplicial maps denoted ${\rm u}^1_2$ : $S'{\cat R} \tdot _S (X) \lra diag {\cat R} \tdotdot \!\!\! _{S',S} (X)$ and ${\rm u}^2_2$ : ${\cat R} \tdot _{S'} (SX) \lra diag {\cat R} \tdotdot \!\!\! _{S',S} (X)$.
The reader should compare the cosimplicial maps ${\rm u}^1_2$, ${\rm u}^2_2$ with the analogous maps $u_2^1$, $u_2^2$ of Appendix \ref{app3}. 

\bigskip

We will prove below that the maps ${\bf R}\lim^\Delta {\rm u}^1_2$ and ${\bf R}\lim^\Delta {\rm u}^2_2$ are isomorphisms in ${\bf ho}{\cat M}$.
Based on that, $\theta_{S', S}$ is defined as the composition of maps and inverses of maps in ${\bf ho}{\cat M}$

\begin{center}
$\xymatrix{
	{\bf R}\lim^\Delta {\cat R}_S \tdot(X) \ar[dd]_{\theta_{S', S}} &
	{\bf R}\lim^\Delta S'{\cat R}_S \tdot(X) \ar[l]_-\cong \ar[d]^-{{\bf R}\lim^\Delta {\rm u}^1_2}_-\cong \\
	& {\bf R}\lim^\Delta diag {\cat R} \tdotdot \!\!\! _{S',S} (X) \\
	{\bf R}\lim^\Delta {\cat R}_{S'} \tdot(X) &
	{\bf R}\lim^\Delta {\cat R}_{S'} \tdot(SX) \ar[l]_-\cong \ar[u]_{{\bf R}\lim^\Delta {\rm u}^2_2}^-\cong &
}$
\end{center}
From this definition of $\theta_{S', S}$, it is easy to see that $\theta_{S', S}$ commutes with the completion maps.

\bigskip

We now prove that the maps ${\bf R}\lim^\Delta {\rm u}^1_2$ and ${\bf R}\lim^\Delta {\rm u}^2_2$ are isomorphisms in ${\bf ho}{\cat M}$.
Using the diagonal argument of Appendix \ref{app2}, note that 
\begin{center}
${\bf R}\lim^\Delta  diag {\cat R} \tdotdot \!\!\! _{S',S} (X) \cong {\bf R} \lim^{\Delta \times \Delta} {\cat R} \tdotdot \!\!\! _{S',S} (X)$
\end{center}

For the bicosimplicial diagram ${\cat R} \tdotdot \!\!\! _{S',S} (X)$ , we will compute ${\bf R} \lim^{\Delta \times \Delta}$ in two ways. 
We denote the first factor of the product category $\Delta \times \Delta$ by $\Delta_{ho}$ (``horizontal'' in diagram \eqref{eqnDoubleComplex}), and the second factor by $\Delta_{ve}$ (``vertical'').

We use first the factorization
\begin{center}
${\bf R}\lim^{\Delta_{ho} \times \Delta_{ve}}$ $\cong$ ${\bf R}\lim^{\Delta_{ho}} {\bf R}\lim^{\Delta_{ve}}$
\end{center}
that computes first ${\bf R}\lim^{\Delta_{ve}}$ in the vertical direction.

In computing ${\bf R}\lim^{\Delta_{ve}} {\cat R} \tdotdot \!\!\! _{S',S} (X)$, we can safely drop the first appearance of $S$ to the left in each of the terms of our diagram (\ref{eqnDoubleComplex}).
After this operation, each column becomes a {\bf contractible} cosimplicial object, by Sec. \ref{secMixTriCot}, (\ref{EqnContractionProperty1}).
We get (by Appendix \ref{app1}, Thm. \ref{thmHomot2}) that ${\bf R}\lim^{\Delta_{ve}} {\cat R} \tdotdot \!\!\! _{S',S} (X) \cong S'{\cat R}_S \tdot(X)$ (as objects of ${\bf ho}({\cat M}^\Delta)$).

\bigskip

It follows that
\begin{center}
${\bf R}\lim^{\Delta_{ho} \times \Delta_{ve}} {\cat R} \tdotdot \!\!\! _{S',S} (X) \cong {\bf R}\lim^{\Delta_{ho}} {\bf R}\lim^{\Delta_{ve}} {\cat R} \tdotdot \!\!\! _{S',S} (X) \cong {\bf R}\lim^\Delta S'{\cat R}_S \tdot(X)$
\end{center}
which proves that the map ${\bf R}\lim^\Delta {\rm u}^1_2$ is an isomorphism in ${\bf ho}{\cat M}$.

\bigskip

Using the factorization in the other order
\begin{center}
${\bf R}\lim^{\Delta_{ho} \times \Delta_{ve}}$ $\cong$ ${\bf R}\lim^{\Delta_{ve}} {\bf R}\lim^{\Delta_{ho}}$
\end{center}
and this time using Sec. \ref{secMixTriCot}, (\ref{EqnContractionProperty2}) in place of Sec. \ref{secMixTriCot}, (\ref{EqnContractionProperty1}) one shows that the other map ${\bf R}\lim^\Delta {\rm u}^2_2$ is an isomorphism in ${\bf ho}{\cat M}$.

\bigskip
\bigskip

So far, we have constructed isomorphisms $\theta_{S', S}$ : $X^{\wedge}_{R,S} \lra X^{\wedge}_{R,S'}$ in ${\bf ho}{\cat M}$.
We are left with proving the cocycle relations $\theta_{S'', S'}\theta_{S', S}$ = $\theta_{S'', S}$ and $\theta_{S, S}$ = $id$.
We will only sketch the proofs of these cocycle relations.

The relation $\theta_{S'', S'}\theta_{S', S}$ = $\theta_{S'', S}$ follows from a triple complex argument - there is nothing essentially new needed for this proof.

\bigskip

The proof for $\theta_{S, S}$ = $id$ is more interesting, as it involves a new idea.
Assume that ($S, \eps, \Delta$) = ($S', \eps, \Delta$).
Using $\Delta$ : $S \lra  S^2$, construct a map of diagrams from the diagram

\bigskip

\begin{equation}
\label{eqnDoubleComplexS}
\xymatrix{ 
	SX \ar[r] \ar[d] &
	SRSX \ar[r] \ar@<1ex>[r] \ar[d] & 
	S(RS)^2X \ar@<1ex>[l] \ar[d] & 
	... \\ 
	SRSX \ar[r] \ar[d] \ar@<1ex>[d] & 
	SRSRSX \ar[r] \ar@<1ex>[r] \ar[d] \ar@<1ex>[d] & 
	SRS(RS)^2X \ar@<1ex>[l] \ar[d] \ar@<1ex>[d] & 
	... \\
	(SR)^2SX \ar[r] \ar[d] \ar@<1ex>[d] \ar@<2ex>[d] \ar@<1ex>[u] & 
	(SR)^2SRSX \ar[r] \ar@<1ex>[r] \ar[d] \ar@<1ex>[d] \ar@<2ex>[d] \ar@<1ex>[u] & 
	(SR)^2S(RS)^2X \ar@<1ex>[l] \ar[d] \ar@<1ex>[d] \ar@<2ex>[d] \ar@<1ex>[u] &
	... \\
	(SR)^3SX \ar[r] \ar@<1ex>[u] \ar@<2ex>[u] & 
	(SR)^3SRSX \ar[r] \ar@<1ex>[r] \ar@<1ex>[u] \ar@<2ex>[u] & 
	(SR)^3S(RS)^2X \ar@<1ex>[l] \ar@<1ex>[u] \ar@<2ex>[u] & 
	... \\
	... & ... & ... & ...
	}
\end{equation}
to the diagram ${\overline {\cat R}}\; \tdotdot \!\!\!\!\! _{S,S} (X)$
\begin{equation}
\label{eqnDoubleComplexS2}
\xymatrix{ 
	S^2X \ar[r] \ar[d] &
	S^2RSX \ar[r] \ar@<1ex>[r] \ar[d] & 
	S^2(RS)^2X \ar@<1ex>[l] \ar[d] & 
	... \\ 
	SRS^2X \ar[r] \ar[d] \ar@<1ex>[d] & 
	SRS^2RSX \ar[r] \ar@<1ex>[r] \ar[d] \ar@<1ex>[d] & 
	SRS^2(RS)^2X \ar@<1ex>[l] \ar[d] \ar@<1ex>[d] & 
	... \\
	(SR)^2S^2X \ar[r] \ar[d] \ar@<1ex>[d] \ar@<2ex>[d] \ar@<1ex>[u] & 
	(SR)^2S^2RSX \ar[r] \ar@<1ex>[r] \ar[d] \ar@<1ex>[d] \ar@<2ex>[d] \ar@<1ex>[u] & 
	(SR)^2S^2(RS)^2X \ar@<1ex>[l] \ar[d] \ar@<1ex>[d] \ar@<2ex>[d] \ar@<1ex>[u] &
	... \\
	(SR)^3S^2X \ar[r] \ar@<1ex>[u] \ar@<2ex>[u] & 
	(SR)^3S^2RSX \ar[r] \ar@<1ex>[r] \ar@<1ex>[u] \ar@<2ex>[u] & 
	(SR)^3S^2(RS)^2X \ar@<1ex>[l] \ar@<1ex>[u] \ar@<2ex>[u] & 
	... \\
	... & ... & ... & ...
	}
\end{equation}

\bigskip

This map from diagram \eqref{eqnDoubleComplexS} to diagram \eqref{eqnDoubleComplexS2} is a pointwise weak equivalence.

Using the formulas of Appendix \ref{app3}, the bicosimplicial object obtained by erasing the top row, and leftmost column in diagram \eqref{eqnDoubleComplexS} is $G_2^* {\cat R} \tdot _S(X)$.
The cosimplicial object on the diagonal of \eqref{eqnDoubleComplexS} is $F_2^* {\cat R} \tdot _S(X)$.

The map from the cosimplicial object on the top row to the cosimplicial object on the diagonal is
\begin{center}
$u_2^1$ : ${\cat R} \tdot _S(X) \lra F_2^* {\cat R} \tdot _S(X)$
\end{center}
and the map from the cosimplicial object on the leftmost column to the diagonal is
\begin{center}
$u_2^2$ : ${\cat R} \tdot _S(X) \lra F_2^* {\cat R} \tdot _S(X)$
\end{center}

By Appendix \ref{app3}, Prop. \ref{proppsikl} it follows that the cosimplicial maps $u_2^1 \simeq u_2^2$ are homotopic, and consequently
\begin{center}
${\bf R}\lim^\Delta {u_2^1}$ = ${\bf R}\lim^\Delta {u_2^2}$ : ${\bf R}\lim^\Delta {\cat R} \tdot _S(X) \lra {\bf R}\lim^\Delta F_2^* {\cat R} \tdot _S(X)$
\end{center}

A diagram chase will quickly show now that $\theta_{S, S}$ = $id$.

\end{proof}

\bigskip

\subsection{The Bousfield-Kan spectral sequence}
\label{subsecBKSpectrSeqCompl}
Throughout Section \ref{subsecBKSpectrSeqCompl}, we assume that ${\cat M}$ is a {\it pointed} model category having a cofibrant-replacement cotriple ($S, \eps, \Delta$).
Consider, as before, a triple ($R, \nu, M$) on ${\cat M}$ with the property that $R$ maps equivalences between cofibrant objects of ${\cat M}$ to equivalences.

We will adapt the homotopy spectral sequence of Chapter \ref{chapHSS}, Sec. \ref{subsecHomotopySpSeq} to construct the Bousfield-Kan spectral sequence.
In what follows, ${\bf Hom}_*$ denotes the pointed function complex defined in Chapter \ref{chapHSS}, Section \ref{secPoiFunCom}.

\bigskip
\bigskip

The Bousfield-Kan spectral sequence is defined as
\begin{center}
$E_r^{s,t}(W, X) = E_r^{s,t}(W, {\cat R}_S \tdot(X))$
$\phantom{foo} (t \geq s \geq 0)$
\end{center}
where on the left we have the homotopy spectral sequence of the cosimplicial object ${\cat R}_S \tdot(X)$ (Chap. \ref{chapHSS}, Sec. \ref{subsecHomotopySpSeq}).

From the description of $E_2$ given by
\begin{center}
$E_2^{s,t}(W, X) = \pi^s \pi_t {\bf Hom}_*(W, {\cat R}_S \tdot(X))$
$\phantom{foo} (t \geq s \geq 0)$
\end{center}
it easily follows that, from $E_2$ on, the Bousfield-Kan spectral sequence is independent of the choice of the cofibrant-replacement cotriple ($S, \eps, \Delta$).

\bigskip
\bigskip

\subsection{Naturality in $R$ of the Bousfield-Kan completion}
\label{subsecNatR}

For Section \ref{subsecNatR}, we return to the general case of a not necessarily pointed model category.

We show that the Bousfield-Kan completion functor exhibits a rigid naturality in the triple ($R, \nu, M$).
Fix the cofibrant-replacement cotriple ($S, \eps, \Delta$) in the model category ${\cat M}$, and assume we have a map $f$ of triples from a triple ($R, \nu, M$) to a triple ($T, \nu, M$) in ${\cat M}$.

We assume that both $R$ and $T$ map equivalences between cofibrant objects to equivalences, so Bousfield-Kan completion can be defined for both of them.

Then we construct in the obvious way a natural map in ${\bf ho}{\cat M}$
\begin{center}
$\xymatrix{ X^{\wedge}_R \ar[r]^{f^\wedge} &
	X^{\wedge}_T
}$
\end{center}

\bigskip

We have the following Theorem and Corollary:
\begin{thm}
\label{thmBKcomplHomotMaps}
If $f$, $g$ are two triple maps from the triple ($R, \nu, M$) to the triple ($T, \nu, M$) then in ${\bf ho}{\cat M}$
\begin{center}
$f^\wedge$ = $g^\wedge$ : $X^{\wedge}_R \lra X^{\wedge}_T$
\end{center}
\end{thm}

\begin{proof}
By Appendix \ref{app1}, Prop. \ref{propMapsOfTriples} there is a natural homotopy of cosimplicial objects
\begin{center}
$f$ $\simeq$ $g$ : ${\cat R}_S \tdot (X) \lra {\cat T}_S \, \!\!\! \tdot (X)$
\end{center}
and by Appendix \ref{app1}, Thm. \ref{thmHomot2} it follows that $f^\wedge$ = $g^\wedge$.
\end{proof}

\begin{cor}
If the triples ($R, \nu, M$), ($T, \nu, M$) admit triple maps in both directions
\begin{center}
$f$ : ($R, \nu, M$) $\lra$ ($T, \nu, M$)

$g$ : ($T, \nu, M$) $\lra$ ($R, \nu, M$)
\end{center}
then there is a canonical isomorphism in ${\bf ho}{\cat M}$ of completions $X^{\wedge}_R \cong X^{\wedge}_T$
\end{cor}

\bigskip
\bigskip
\bigskip

                % Local Variables: 
                % mode: latex 
                % tex-command: amslatex-command
                % tex-dvi-view-command: tex-oneside-view-command
                % tex-start-of-header: "\\special"
                % tex-end-of-header: "\\maketitle"
                % abbrev-mode: t
                % End:

%% file: ch1.3.tex
\section{An Example: Completion of Simplicial Algebras}
\label{secExamplesBKCompletion}
\bigskip
\bigskip
\bigskip

In this Section we describe one example of Bousfield-Kan completion: the completion of simplicial commutative algebras, with respect to certain triples that include the abelianization triple.

In order to do that, we first recall the model category structure on the category ${\cat A}_\alpha$ of simplicial commutative algebras $X.$ of the form
\begin{center}
\[ \xymatrix {
	A \ar[r] \ar@/_1pc/[rr]_{\alpha}
	& X. \ar[r]
	& B
	}
\]
\end{center} 
where $\alpha$ is a fixed commutative algebra map.

Then, for a fixed commutative algebra map $\beta$ : $B \lra C$, we consider a certain adjoint pair between the category ${\cat A}_\alpha$ and the category $sC$-$mod$ of simplicial $C$-modules (see Subsection \ref{secQuiAdjPair}). 
This adjoint pair gives rise to a triple ($R, \nu, M$) on the category of simplicial commutative algebras ${\cat A}_\alpha$.
For a simplicial algebra $X.$ the simplicial algebra $RX.$ is going to be the algebra extension with $B$ of the $C$-module $\Omega_{X.|A} \tens_A C$, namely $RX.$ = $B \bigoplus (\Omega_{X.|A} \tens_A C)$, where $\Omega$ denotes the K\"{a}hler differentials.

In particular, if $B$ = $C$, the functor $R$: ${\cat A}_\alpha \lra {\cat A}_\alpha$ that we construct is just {\it the abelianization functor}, seen as a functor from the category ${\cat A}_\alpha$ to itself.

We define in our context the Bousfield-Kan completion with respect to the triple ($R, \nu, M$), as in Section \ref{secBCCF}.

Our conjecture is that simplicial augmented connected $k$-algebras ($k$ a field) are complete with respect to abelianization, and their associated absolute Bousfield-Kan spectral sequence is convergent. 

\bigskip
\bigskip

\subsection{${\cat A}_\alpha$ as a model category}
${\cat A}_\alpha$ carries a cofibrantly generated model category structure, constructed as follows.
We have an adjoint pair of functors 
\begin{center}
\[ \xymatrix {
	sSets \ar@<.5ex>[rr]^-{free}
	& & {\cat A}_\alpha \ar@<.5ex>[ll]^-{forget}
	}
\]
\end{center} 

$sSets$ is a cofibrantly generated model category.
Using the above adjoint pair of functors ($free$, $forget$) we can lift (in the sense of \cite{Dwyer-Kan-Hirschhorn}, 9.1 and 9.9) the model category structure on $sSets$ to a cofibrantly generated model category structure on ${\cat A}_\alpha$, by saying that a simplicial algebra map
\begin{center}
$X. \ra Y.$
\end{center} 
is a weak equivalence (resp. a fibration) if its underlying map of simplicial sets is a weak equivalence (resp. a fibration).

\bigskip
\bigskip

\subsection{The Quillen adjoint pair between ${\cat A}_\alpha$ and $sC$-$mod$}
\label{secQuiAdjPair}
Let us fix a commutative algebra map $\beta$ : $B \lra C$.

\bigskip

If $M$ is a $C$-module ($M \in Ob(C$-$mod)$), we construct a $B$-algebra structure on the $B$-module direct sum
\begin{center}
$B \bigoplus M$
\end{center} 
by defining multiplication as $(b_1,x_1) \cdot (b_2,x_2) = (b_1 b_2, b_1 x_2 + b_2 x_1)$. 
$B \bigoplus M$ is an augmented $A$-algebra in the obvious way.

This construction, extended in the obvious way for simplicial $C$-modules, defines a functor classically denoted $+$ : $sC$-$mod \lra {\cat A}_\alpha$, from the category of simplicial $C$-modules to the category ${\cat A}_\alpha$.

\bigskip

In the particular case $B$ = $C$, it is easy to see that the functor $+$ is an embedding, and the image of $+$ is equivalent to $Ab{\cat A}_\alpha$, the category of abelian objects of ${\cat A}_\alpha$. 

In general, the functor $+$ admits a left adjoint denoted $Q_C$ (if $B$ = $C$, $Q_C$ is just the {\bf abelianization functor} $Ab$):
\begin{center}
\bigskip 
${\cat A}_\alpha
	\xymatrix {
	 \ar@<.5ex>[r]^-{Q_C}
	& \ar@<.5ex>[l]^-{+}
	}
sC {\mathrm -}mod$

\bigskip 

$Q_C X.$ = $\Omega_{X.|A} \tens_A C$

\bigskip 

$M_+$ = $B \bigoplus M.$

\bigskip 
\end{center} 
where $\Omega$ denotes the K\"{a}hler differentials.

The category $sC$-$mod$ is also a model category (again by \cite{Dwyer-Kan-Hirschhorn}, 9.1 and 9.9), with the property that a simplicial $A$-module map
\begin{center}
$M. \ra N.$
\end{center} 
is a weak equivalence (resp. a fibration) if its underlying map of simplicial sets is a weak equivalence (resp. a fibration).

In fact, all objects of $sC$-$mod$ have an underlying structure of simplicial (abelian) groups, therefore their underlying simplicial set is fibrant.
Consequently, all objects of $sC$-$mod$ are fibrant.

With these model category structures on ${\cat A}_\alpha$ and $sC$-$mod$, the adjoint pair of functors ($Q_C$, $+$) is a Quillen adjoint pair.

\bigskip

\subsection{The completion functor for simplicial algebras}

Denote $R$ the composition
\begin{center}
$R = + \comp Q_C$ : ${\cat A}_\alpha \lra {\cat A}_\alpha$.
\end{center} 
Since ($Q_C$, $+$) is an adjoint pair of functors, $R$ is part of a triple ($R, \nu, M$) in a natural way.

\bigskip

We claim that all hypotheses (Section \ref{subsecConstOfCompl}) necessary to construct the Bousfield-Kan completion with respect to $R$ are satisfied.

Indeed, ($Q_C$, $+$) is a Quillen adjoint pair, and all objects of $sC$-$mod$ are fibrant: it follows easily that $R$ carries equivalences between cofibrant objects to equivalences.
Furthermore, ${\cat A}_\alpha$ is a cofibrantly generated model category, so it has a cofibrant-replacement cotriple ($S, \eps, \Delta$).

\bigskip

If $X. \in Ob{\cat A}_\alpha$ is a simplicial algebra, its Bousfield-Kan completion with respect to the triple ($R, \nu, M$) is

\begin{center}
$X^{\wedge}_R = {\mathbf R}\lim^{\cat \Delta} {\cat R}_S \tdot(X.)$
\end{center}

If $B$ = $C$, the functor $R$ is just the abelianization functor regarded as a functor ${\cat A}_\alpha \lra {\cat A}_\alpha$, and in this case we denote the completion as $X^{\wedge}_{Ab}$ -in effect, it is the Bousfield-Kan completion with respect to abelianization.

\bigskip

If $A$ = $B$, the model category ${\cat A}_{id_A}$ is {\it pointed} and, given a second simplicial algebra $W. \in Ob{\cat A}_{id_A}$, we have the Bousfield-Kan spectral sequence $E_r^{s,t}(W., X.)$, as described in Section \ref{subsecBKSpectrSeqCompl}.
This spectral sequence computes the homotopy groups of the function complex ${\bf Hom}_*(W., X^{\wedge}_R)$.

\bigskip

\subsection{A conjecture about connected simplicial algebras}

Suppose that in the constructions above we take $A$ = $B$ = $C$ = $k$ a field.
In this case, we denote by ${\cat A}_k$ the category ${\cat A}_{id_k}$.

The natural map $X. \lra X^{\wedge}_{Ab}$ is an equivalence when $X.$ is an abelian object in ${\cat A}_k$, since ${\overline {\cat R}}\; \tdot \!\!\! _S(X.)$ is contractible for $X.$ abelian.

By definition, we say that $X.$ is a {\bf connected} simplicial algebra if the natural map
\begin{center}
$k \ra \pi_0 (X.)$ 
\end{center}
is an isomorphism.
 
We state the following
\begin{conj}
If $X.$ is a connected simplicial commutative augmented $k$-algebra, then the natural map
\begin{center}
$X. \lra X^{\wedge}_{Ab}$ 
\end{center}
is an equivalence, and the {\bf absolute} Bousfield-Kan spectral sequence $E_r^{s,t}(k[T], X.)$ is convergent.
\end{conj}

$k[T]$ is just the polynomial ring in one variable, and should be thought of as the analogue of the $0$-sphere in ${\cat A}_k$.

The Bousfield-Kan spectral sequence $E_r^{s,t}(k[T], X.)$ can be equivalently constructed as follows.
Take ${\cat R}_S \tdot(X.)$ and regard it as a cosimplicial space, just by forgetting the algebra structure.
$E_r^{s,t}(k[T], X.)$ is the absolute homotopy spectral sequence of the cosimplicial space ${\cat R}_S \tdot(X.)$: it is a truncation of the spectral sequence of the bi-complex of abelian groups underlying ${\cat R}_S \tdot(X.)$.

\bigskip
\bigskip
\bigskip

                % Local Variables: 
                % mode: latex 
                % tex-command: amslatex-command
                % tex-dvi-view-command: tex-oneside-view-command
                % tex-start-of-header: "\\special"
                % tex-end-of-header: "\\maketitle"
                % abbrev-mode: t
                % End:

%% file: chap2.tex
\chapter{The Homotopy Spectral Sequence for a Model Category}
\label{chapHSS}

\bigskip
\bigskip
\bigskip

This Chapter is a reference for the general model category results used in the rest of these notes.
The end result of this Chapter is the construction of a homotopy spectral sequence of a cosimplicial object in a pointed model category, as a generalization of the homotopy spectral sequence of a cosimplicial space of Bousfield and Kan, \cite{Bousfield-Kan1}.

\bigskip

In Section \ref{secCoreThms}, we recall the basics behind homotopy limits and the function complex.
Most important, the function complex has the property that it commutes with homotopy limits.

\bigskip

The homotopy spectral sequence of a cosimplicial object $X \tdot$ in a pointed model category is defined (in Section \ref{secSpSeq}) just as the classical homotopy spectral sequence of the cosimplicial pointed space ${\bf Hom}_*(W, X \tdot)$.
In this formula, $W$ is an arbitrary object in the model category, and ${\bf Hom}_*($-,-$)$ denotes the {\it pointed} function complex associated to the pointed model category. 

So, the homotopy spectral sequence is actually the spectral sequence of a certain $tot$ tower of pointed spaces constructed using the pointed function complex.
We will show that the $tot$ tower in discussion in fact comes from a tower in the model category.

To that end, in Section \ref{sectottot^n}, we recall the basic properties of $tot$ and $tot^n$, and define for a model category functors ${\bf tot}$ and ${\bf tot}^n$, and dually  ${\bf cotot}$ and ${\bf cotot}_n$.
There is going to be an essential difference between $tot$, $tot^n$ and ${\bf tot}$, ${\bf tot}^n$, ${\bf cotot}$ and ${\bf cotot}_n$: the first ones yield objects in the model category $sSets$, whereas the last ones yield objects in the homotopy category ${\bf ho}{\cat M}$ of a model category.

\bigskip

Finally, Section \ref{secSpSeq} deals with the homotopy spectral sequence itself.

\bigskip
\bigskip
\bigskip

\section{Homotopy Limits. The Function Complex}
\label{secCoreThms}
\bigskip
\bigskip

In this Section we collect general results about model categories, that we will need later in constructing the homotopy spectral sequence.
The material we present is collated from the excellent Dwyer, Kan and Hirschhorn monograph \cite{Dwyer-Kan-Hirschhorn}, and from the Dwyer, Kan articles \cite{Dwyer-Kan1}, \cite{Dwyer-Kan2}, \cite{Dwyer-Kan3}.

In Section \ref{subsecHomotLim} we define homotopy limits as derived functors of limits.
For properties of homotopy limits, we direct the reader to the monographs \cite{Dwyer-Kan-Hirschhorn}, \cite{Hirschhorn}.

In Section \ref{subsecTheFuncComplex} we define the function complex, for a pair of small categories (${\cat C}$, ${\cat W}$) with $Ob{\cat C}$ = $Ob{\cat W}$.
If ${\cat C}$ = ${\cat M}$ is a (not necessarily small) model category, with ${\cat W}$ the subcategory of weak equivalences, then the function complex is ``homotopically small'', in a sense made precise in that Section.

In Section \ref{subsecHomLimFunComples}, we observe that for a model category, the function complex commutes with homotopy limits (this is a result due to Dwyer, Kan and Hirschhorn).

\bigskip

In order to define the homotopy spectral sequence for a pointed model category, we need a pointed version of the function complex.
In Section \ref{secPoiFunCom}, we develop the {\it pointed} function complex, defined for a pair of categories (${\cat C}$, ${\cat W}$) such that ${\cat C}$ is pointed.

Last, in Section \ref{subsecHtpyLimPoiFuncComplex} we observe that, for a pointed model category, the pointed function complex commutes with pointed homotopy limits.

\bigskip
\bigskip

Let ${\cat M}$ be a model category.
Let ${\cat W}$ be the subcategory of weak equivalences of ${\cat M}$.
${\cat W}$ has the same objects as ${\cat M}$, and has as maps the weak equivalences of ${\cat M}$.

\subsection{Homotopy Limits}
\label{subsecHomotLim}

The homotopy category ${\bf ho}{\cat M}$ of ${\cat M}$ is defined as the localization ${\cat M}[{\cat W}^{-1}]$ of ${\cat M}$ with respect to ${\cat W}$.

If ${\cat D}$ is a small category, we denote by ${\cat M}^{\cat D}$ the category of ${\cat D}$-diagrams in ${\cat M}$.
We define weak equivalences in ${\cat M}^{\cat D}$ as ${\cat W}^{\cat D}$.
Denote ${\bf ho}({\cat M}^{\cat D})$=${\cat M}^{\cat D}[({\cat W}^{\cat D})^{-1}]$.

One proves that ${\bf ho}{\cat M}$, as well as ${\bf ho}({\cat M}^{\cat D})$, are categories within the universe we start with (\cite{Dwyer-Kan-Hirschhorn}).

\bigskip

Since ${\cat M}$ is complete and cocomplete, we have adjoint pairs of functors:
\begin{center}
\bigskip
${\colim}^{\cat D} : {\cat M}^{\cat D} \xymatrix{ \ar@<.5ex>[r] & \ar@<.5ex>[l] } {\cat M} : c$

\bigskip
$c : {\cat M} \xymatrix{ \ar@<.5ex>[r] & \ar@<.5ex>[l] } {\cat M}^{\cat D} : {\lim}^{\cat D}$
\bigskip
\end{center}
where the functor $c$ forms the constant diagram.

One shows (see \cite{Dwyer-Kan-Hirschhorn}) that the functors ${\colim}$ and ${\lim}$ admit a left, respectively right derived functor with respect to localization, that still form adjoint pairs.

\bigskip
\begin{center}
${\bf L}{\colim}^{\cat D} : {\bf ho}({\cat M}^{\cat D}) \xymatrix{ \ar@<.5ex>[r] & \ar@<.5ex>[l] } {\bf ho}{\cat M} : c$

\bigskip
$c : {\bf ho}{\cat M} \xymatrix{ \ar@<.5ex>[r] & \ar@<.5ex>[l] } {\bf ho}({\cat M}^{\cat D}) : {\bf R}{\lim}^{\cat D}$
\end{center}
\bigskip

Please observe that if a functor $F$ : ${\cat M} \lra {\cat N}$ carries weak equivalences to weak equivalences then the diagram
\begin{center}
\bigskip
$\xymatrix{
	{\cat M} \ar[r]^F \ar[d] &
	{\cat N} \ar[d] \\
	{\bf ho}{\cat M} \ar@<.5ex>[r]^{{\bf L}F} \ar@<-.5ex>[r]_{{\bf R}F} &
	{\bf ho}{\cat N}
}$
\bigskip
\end{center}
commutes.
In this situation, we will denote again by $F$=${\bf L}F$=${\bf R}F$ the localization of $F$: we hope the context will make it clear if we refer to the localization of $F$ or to $F$ itself.

\bigskip

It is not true that ${\colim}$ of a weak equivalence between two pointwise cofibrant diagrams is a weak equivalence.
This would be a most desirable property.
In \cite{Dwyer-Kan-Hirschhorn} a different functor $hocolim^{\cat D} : {\cat M}^{\cat D} \ra {\cat M}$ is constructed, using the choice of a cosimplicial framing. 
$hocolim$ has this ``desirable'' property, and furthermore it satisfies ${\bf L}hocolim^{\cat D} $ $\cong$ $ {\bf L}{\colim}^{\cat D}$.

Dually, choosing a simplicial framing one constructs $holim^{\cat D} : {\cat M}^{\cat D} \ra {\cat M}$ that carries weak equivalences between pointwise fibrant diagrams to weak equivalences, with the property that ${\bf R}holim^{\cat D} $ $\cong$ $ {\bf R}{\lim}^{\cat D}$.

\bigskip

\subsection{The function complex}
\label{subsecTheFuncComplex}

Consider a pair of small categories (${\cat C}$, ${\cat W}$) (this means, ${\cat W}$ a subcategory of ${\cat C}$) with $Ob{\cat C}=Ob{\cat W}$.
When ${\cat C} = {\cat M}$ is a model category, we always take ${\cat W}$ to be the subcategory of weak equivalences.

We will construct the function complex ${\bf Hom}(X,Y)$ between any two objects $X$, $Y$ of ${\cat C}$ as a homotopy type of a simplicial set.

To that end, we outline Dwyer and Kan's construction of the hammock localization $L^H{\cat C}$ (see \cite{Dwyer-Kan2}).
Even though it is not apparent from the notation, the hammock localization depends both on ${\cat C}$ and ${\cat W}$.

\bigskip

Denote ${\cat O}$ the set $Ob{\cat C}=Ob{\cat W}$ (it is a set, since ${\cat C}$ and ${\cat W}$ are small categories) .

The hammock localization $L^H{\cat C}$ is a category enriched over simplicial sets, with $Ob (L^H{\cat C})$ = ${\cat O}$.
It is defined as follows:
The $k$-simplices of the simplicial set $Hom_{L^H{\cat C}}(X, Y)$ are the commutative diagrams in ${\cat C}$ of the form

\begin{center}
\[
\xymatrix{
	& C_{0,1} \ar[d]^\sim \ar@{-}[r]
	& C_{0,2} \ar[d]^\sim \ar@{-}[r]
	& ... \ar@{-}[r]
	& C_{0,n-1} \ar[d]_\sim \ar@{-}[dr] 
	& \\
	X \ar@{-}[ur] \ar@{-}[r] \ar@{-}[ddr] 
	& C_{1,1} \ar[d]^\sim \ar@{-}[r]
	& C_{1,2} \ar[d]^\sim \ar@{-}[r]
	& ... \ar@{-}[r]
	& C_{1,n-1} \ar[d]_\sim \ar@{-}[r] 
	& Y \\
	& ... \ar[d]^\sim & ... \ar[d]^\sim & & ... \ar[d]_\sim & \\
	& C_{k,1} \ar@{-}[r]
	& C_{k,2} \ar@{-}[r]
	& ... \ar@{-}[r]
	& C_{k,n-1} \ar@{-}[uur] 
	& 
} \]
\end{center}

\begin{list}{-}{where}
  \item $n$ is any integer $\ge 0$
  \item the vertical maps are in ${\cat W}$  
  \item in each column, all maps go in the same direction; if they go left, they are in ${\cat W}$
  \item the maps in adjacent columns go in opposite directions
  \item no column contains only identity maps.
\end{list}

\bigskip

We have a simplicial functor from the category ${\cat C}$ to the hammock localization $L^H({\cat C}, {\cat W})$  
\begin{center}
$F$: ${\cat C} \lra L^H{\cat C}$
\end{center}
defined by $F(X)$ = $X$, by $F(f:X \ra Y)$ = the diagram $X \xymatrix{\ar[r]^f &} Y$ and by $F(id_X)$ = ``the diagram'' $X$.

\bigskip

We define the function complex ${\bf Hom}(X,Y)$ as the {\it homotopy type} of the simplicial set $Hom_{L^H{\cat C}}(X,Y)$.

Using the simplicial functor $F$, we will regard ${\bf Hom}($-$,$-$)$ as a homotopy type of diagrams of simplicial sets in $sSets^{{\cat C}^{op}\times{\cat C}}$, that is, an object of ${\bf ho}(sSets^{{\cat C}^{op}\times{\cat C}})$.

\bigskip

\bigskip

If ${\cat C} = {\cat M}$ is a model category, $Ob{\cat M}=Ob{\cat W}$ is a class and of course there's no problem in performing these constructions in the category $Classes$ of proper classes (a category in a higher universe) as opposed to the category $Sets$.

We have:
\begin{thm}
\label{thmC1}
For an arbitrary model category ${\cat M}$, fix a framing. 
Then there exist natural isomorphisms in ${\bf ho}(sClasses^{{\cat M}^{op}\times{\cat M}})$
\begin{center}
${\bf Hom}(X,Y) $ $\cong$ $ Hom_{\cat M}({\bf X \tdot}, Y') $ $\cong$ $
Hom_{\cat M}(X', {\bf Y .}) $ $\cong$ $ {\bf diag} Hom_{\cat M}({\bf X \tdot}, {\bf Y .})$
\end{center}

Here $X'$ is a canonical cofibrant replacement of $X$, ${\bf X \tdot}$ the cosimplicial frame of $X'$, $Y'$ is a canonical fibrant replacement of $Y$ and ${\bf Y .}$ the simplicial frame of $Y'$.
${\bf diag}$ denotes the diagonal of a bisimplicial set.
\end{thm}

\begin{proof}
See \cite{Dwyer-Kan3}, 4.4. and its proof.
\end{proof}

As a consequence, for a not necessarily small model category ${\cat M}$, ${\bf Hom}(X,Y)$ is {\bf homotopically small} for $X$, $Y \, \eps \, Ob({\cat M})$, that is, it has the homotopy type of a simplicial {\bf set}.

Even more than that is true.
If ${\cat M}$ is a model category, by Theorem \ref{thmC1}, the function complex ${\bf Hom}($-$,$-$) \in Ob({\bf ho}(sClasses^{{\cat M}^{op}\times{\cat M}}))$ descends to an object in ${\bf ho}(sSets^{{\cat M}^{op}\times{\cat M}})$.

This descent is unique in the following weak sense.
If ${\cat D}$, ${\cat D}'$ are small categories with functors ${\cat D} \lra {\cat M}$ and ${\cat D'} \lra {\cat M}$, then the restriction of the function complex ${\bf Hom}($-$,$-$)$ to ${\bf ho}(sClasses^{{\cat D}^{op}\times{\cat D}'})$ descends in a unique way to ${\bf ho}(sSets^{{\cat D}^{op}\times{\cat D}'})$.

From now on, we will always consider the function complex ${\bf Hom}($-$,$-$)$ associated to a model category ${\cat M}$ as an object of ${\bf ho}(sSets^{{\cat M}^{op}\times{\cat M}})$.

\bigskip

\subsection{Homotopy limits commute with the function complex}
\label{subsecHomLimFunComples}
The following theorem, due to Dwyer, Kan and Hirschhorn \cite{Dwyer-Kan-Hirschhorn}, unravels most of the homotopical properties of model categories: 

\bigskip

\begin{thm}
\label{thmC2}
Let ${\cat D}$, ${\cat D}'$ be small categories, $X \tdot$ an object of ${\cat M}^{\cat D}$ and $Y \tdot$ an object of ${\cat M}^{{\cat D}'}$.
Then the natural map below is an isomorphism in ${\bf ho}(sSets)$:
\begin{equation}
\label{eqnFuncComplHomotLim}
{\bf Hom}({\bf L}{\colim}^{\cat D}  X \tdot, {\bf R}{\lim}^{{\cat D}'} Y \tdot) \xymatrix{\ar[r]^\cong &} {\bf R}{\lim}^{{\cat D}^{op} \times {\cat D}'} {\bf Hom}(X \tdot, Y \tdot)
\end{equation}
\end{thm} 

\bigskip

\begin{note}
${\bf Hom}$ carries weak equivalences in ${\cat M}^{op} \times {\cat M}$ to weak equivalences of simplicial sets.
Observe that, using our convention, on the left hand side of (\ref{eqnFuncComplHomotLim}), ${\bf Hom}$ stands for the localization of the functor denoted also ${\bf Hom}$.
\end{note}
\begin{note}
On the right hand side, the ${\cat D}^{op} \times {\cat D}'$-diagram ${\bf Hom}(X \tdot, Y \tdot)$ is defined only up to weak equivalence in $sSets^{{\cat D}^{op}\times{\cat D}'}$.
The natural map (\ref{eqnFuncComplHomotLim}) is constructed using the adjunction
\begin{center}
$c : {\bf ho}sSets \xymatrix{ \ar@<.5ex>[r] & \ar@<.5ex>[l] } {\bf ho}(sSets^{{\cat D}^{op} \times {\cat D}'}) : {\bf R}{\lim}^{{\cat D}^{op} \times {\cat D}'}$
\end{center}
applied to the natural map in ${\bf ho} Sets^{{\cat D}^{op}\times{\cat D}'}$
\begin{center}
$c{\bf Hom}({\bf L}{\colim}^{\cat D}  X \tdot, {\bf R}{\lim}^{{\cat D}'} Y \tdot) \lra {\bf Hom}(X \tdot, Y \tdot)$
\end{center}
\end{note}
\begin{note}
The map (\ref{eqnFuncComplHomotLim}) is natural with respect to functors $F$: ${\cat D}_1 \ra {\cat D}$ and $F'$: ${{\cat D}'\!\!}_1 \ra {\cat D}'$, meaning that it makes the following diagram commutative:

\begin{center}
$\xymatrix { 
	{{\bf Hom}({\bf L}{\colim}^{\cat D}  X \tdot, {\bf R}{\lim}^{{\cat D}'} Y \tdot)} \ar[r]^{{\cong}} \ar[d] &
	{{\bf R}{\lim}^{{\cat D}^{op} \times {\cat D}'} {\bf Hom}(X \tdot, Y \tdot)} \ar[d] \\
	{{\bf Hom}({\bf L}{\colim}^{{\cat D}_1}  F^{*}X \tdot, {\bf R}{\lim}^{{{\cat D}'\!\!}_1} {F'}^{*}Y \tdot)} \ar[r]^{{\cong}} &
	{{\bf R}{\lim}^{{\cat D}^{op}_1 \times {{\cat D}'\!\!}_1} {\bf Hom}(F^{*}X \tdot, {F'}^{*}Y \tdot)} 
}$
\end{center}
We refer to this property as ``naturality in ${\cat D}$ and ${\cat D}'$''.
\end{note}

\bigskip

\begin{proof}[Proof of Theorem \ref{thmC2}]
The proof can be split into two steps.

Step 1: The proof reduces to the case when either of ${\cat D}$ and ${\cat D}'$ is a singleton category.
If we assume the result known for ${\cat D}$ or ${\cat D}'$ singleton categories, then

\begin{center}
${\bf Hom}({\bf L}{\colim}^{\cat D}  X \tdot, {\bf R}{\lim}^{{\cat D}'} Y \tdot) \cong {\bf R}{\lim}^{{\cat D}^{op}} {\bf Hom}(X \tdot, {\bf R}{\lim}^{{\cat D}'} Y \tdot) \cong$

${\bf R}{\lim}^{{\cat D}^{op}} {\bf R}{\lim}^{{\cat D}'} {\bf Hom}(X \tdot, Y \tdot) \cong {\bf R}{\lim}^{{\cat D}^{op} \times {\cat D}'} {\bf Hom}(X \tdot, Y \tdot)$
\end{center}
(last isomorphism because ${\bf R}{\lim}^{{\cat D}^{op}} {\bf R}{\lim}^{{\cat D}'} $ $\cong$ $ {\bf R}{\lim}^{{\cat D}^{op} \times {\cat D}'}$ , as a consequence of \cite{Dwyer-Kan-Hirschhorn}, 58.3 and 56.5).

\bigskip

Step 2: Let's assume ${\cat D}'$ a singleton category (the other case is similar).
We want to show that there is a natural isomorphism in ${\bf ho}{\cat M}$

\begin{center}
$ {\bf Hom}({\bf L}{\colim}^{\cat D}  X \tdot, Y) \xymatrix{ \ar[r]^\cong &} {\bf R}{\lim}^{{\cat D}^{op}} {\bf Hom}(X \tdot, Y) $
\end{center}

This is proved by Dwyer, Kan and Hirschhorn \cite{Dwyer-Kan-Hirschhorn}, 62.2, using particular representations of the function complex and of homotopy limits.
We fix a cosimplicial framing in ${\cat M}$, and chooose a representative of ${\bf Hom}($-$,$-$)$ in $sSets^{{\cat M}^{op}\times{\cat M}}$, denoted $Hom_\ell($-$,$-$) \in Ob(sSets^{{\cat M}^{op}\times{\cat M}})$ as in Theorem \ref{thmC1}: 
\begin{center}
$
Hom_\ell(A, B) = Hom_{\cat M}({\bf A \tdot}, B')
$
\end{center}
Here $A$, $B \in Ob{\cat M}$, ${\bf A \tdot}$ is the cosimplicial frame of a canonical cofibrant replacement of $A$, and $B'$ is a canonical fibrant replacement of $B$.

Dwyer, Kan and Hirschhorn prove that there is a natural weak equivalence
\begin{equation}
\label{eqn4thmC2}
Hom_\ell (hocolim^{i \in {\cat D}}  (X^i)', Y') \simeq holim^{i \in {\cat D}^{op}} Hom_\ell ((X^i)', Y')
\end{equation}
where $(X^i)'$ is a canonical cofibrant replacement of $X^i$, and $Y'$ a canonical fibrant replacement of $Y$.

The weak equivalence \eqref{eqn4thmC2}, passed to the homotopy category ${\bf ho}{\cat M}$, is actually the isomorphism we are looking for:

The left hand side of \eqref{eqn4thmC2} is ${\bf Hom}({\bf L}{\colim}^{\cat D}  X \tdot, Y)$.

The simplicial set $Hom_\ell ((X^i)', Y')$ is fibrant because e.g. of \cite{Dwyer-Kan-Hirschhorn}, 52.3.  
Therefore the right hand side of \eqref{eqn4thmC2} is ${\bf R}{\lim}^{{\cat D}^{op}} {\bf Hom}(X \tdot, Y)$.
\end{proof}

\bigskip

\subsection{The pointed function complex}
\label{secPoiFunCom}

Consider a pair of small categories (${\cat C}$, ${\cat W}$) with $Ob{\cat C}=Ob{\cat W}$, with the property that ${\cat C}$ is {\it pointed}.

The simplicial functor from ${\cat C}$ to the hammock localization 
\begin{center}
$F$: ${\cat C} \lra L^H{\cat C}$
\end{center}
defines basepoints in the simplicial hom-sets of $L^H{\cat C}$, so the hammock localization becomes a category enriched over {\it pointed} simplicial sets.
We denote $sSets_*$ the category of pointed simplicial sets.

Taking account of the basepoints, we define the {\it pointed} function complex of the pair (${\cat C}$, ${\cat W}$) as a homotopy type of a ${\cat M}^{op}\times{\cat M}$ diagram of pointed simplicial sets
\begin{center}
${\bf Hom}_*(-,-) \in Ob ({\bf ho}(sSets_*^{{\cat M}^{op}\times{\cat M}}))$ 

${\bf Hom}_*(X, Y) = Hom_{L^H{\cat C}}(X, Y)$
\end{center}

We have the following extension to Theorem \ref{thmC1}:
\begin{thm}
\label{thmC1.1}
For a pointed framed model category ${\cat M}$, the natural map below is an isomorphisms in ${\bf ho}(sClasses_*^{{\cat M}^{op}\times{\cat M}})$
\begin{center}
${\bf Hom}_*(X,Y) $ $\cong$ $ Hom_{\cat M}({\bf X \tdot}, Y')$ $\cong$
$Hom_{\cat M}(X', {\bf Y .}) $ $\cong$ $ {\bf diag} Hom_{\cat M}({\bf X \tdot}, {\bf Y .})$
\end{center}

The basepoints in the last three terms are given by zero maps.

Here $X'$ is a canonical cofibrant replacement of $X$, ${\bf X \tdot}$ the cosimplicial frame of $X'$, $Y'$ is a canonical fibrant replacement of $Y$ and ${\bf Y .}$ the simplicial frame of $Y'$.

\end{thm}

\begin{proof}
One just keeps track of basepoints in the proof of Theorem \ref{thmC1}.
\end{proof}

As before, if ${\cat M}$ is a pointed model category, by Theorem \ref{thmC1.1}, the pointed function complex ${\bf Hom}_*($-$,$-$) \in Ob({\bf ho}(sClasses_*^{{\cat M}^{op}\times{\cat M}}))$ descends (uniquely, in the same weak sense we explained before) to an object in ${\bf ho}(sSets_*^{{\cat M}^{op}\times{\cat M}})$.
From now on, we will always consider the pointed function complex ${\bf Hom}_*($-$,$-$)$ associated to a model category ${\cat M}$ as an object of ${\bf ho}(sSets_*^{{\cat M}^{op}\times{\cat M}})$.

\subsection{Homotopy limits commute with the pointed function complex}
\label{subsecHtpyLimPoiFuncComplex}
Theorem \ref{thmC2} extends to

\begin{thm}
\label{thmC2.1}
Let ${\cat M}$ be a pointed model category.

Let ${\cat D}$, ${\cat D}'$ be small categories, $X \tdot$ an object of ${\cat M}^{\cat D}$ and $Y \tdot$ an object of ${\cat M}^{{\cat D}'}$.
Then the natural map below is an isomorphism in ${\bf ho}(sSets_*)$:
\begin{center}
${\bf Hom}_*({\bf L}{\colim}^{\cat D}  X \tdot, {\bf R}{\lim}^{{\cat D}'} Y \tdot) \xymatrix{\ar[r]^\cong &} {\bf R}{\lim}^{{\cat D}^{op} \times {\cat D}'} {\bf Hom}_*(X \tdot, Y \tdot)$
\end{center}
\end{thm} 

\begin{proof}
Follows by keeping track of basepoints in the proof of Theorem \ref{thmC2}.
The underlying space of ${\bf R}{\lim}$ of pointed spaces is (weakly equivalent to) ${\bf R}{\lim}$ of the underlying spaces. 
\end{proof}

\bigskip
\bigskip
\bigskip

\section{$tot$, $tot^n$ . Construction as Homotopy Limits}
\label{sectottot^n}
\bigskip
\bigskip

In this Section, we recover classical results about $tot$ and $tot^n$ of simplicial sets, as seen through the prism of Section \ref{secCoreThms}.

The material in this Section is preparatory for the description of the homotopy spectral sequence in Section \ref{secSpSeq}.
\bigskip

We define $tot$ and $tot^n$ in Section \ref{subsectottotn}.

$tot^n X \tdot$ is better understood as $tot$ of the $n$-th coskeleton of $X \tdot$. 
In Section \ref{subsecCosk} we speak of skeleton and coskeleton of cosimplicial objects in a (complete and cocomplete) category, and then specialize to the case of cosimplicial spaces to deduce properties of $tot$ and $tot^n$.

In Section  \ref{subsecTotTotnHomotProp}, we focus on homotopical properties of $tot$ and $tot^n$.
Among others, for a Reedy fibrant cosimplicial space $X \tdot$ we show that $tot X \tdot$ computes the homotopy limit of $X \tdot$ seen as a cosimplicial diagram, and $tot^n X \tdot$ computes the homotopy limit of the $n$-truncation of the same cosimplicial diagram $X \tdot$.

Section \ref{subsecpointedtotcotot} deals with the pointed versions of $tot$ and $tot^n$, denoted $tot_*$ and $tot^n_*$.

\bigskip

In Section \ref{subsecTotCotot}, the focus shifts to general model categories.
We define functors 
\begin{center}
${\bf tot}$, ${\bf tot}^n$ : ${\bf ho}({\cat M}^\Delta) \lra {\bf ho}{\cat M}$
\end{center}
as generalizations for model categories of {\bf the right derived functors} of $tot$ and $tot^n$.
These functors are later needed for describing the homotopy spectral sequence of a cosimplicial object in a pointed model category.
By duality, we also define functors ${\bf cotot}$ and ${\bf cotot}_n$.

\bigskip
\bigskip

\subsection{$tot$ and $tot^n$}
\label{subsectottotn}
For $X \tdot , Y \tdot \in Ob(sSets^\Delta)$ cosimplicial spaces, $hom(X \tdot,Y \tdot)$ denotes the simplicial hom-set given by 

\begin{center}
$hom(X \tdot,Y \tdot)_k$ = cosimplicial maps $X \tdot \times \Delta [k] \ra Y \tdot$ 
\end{center}

We denote $\Delta \tdot$ the cosimplicial standard space (Bousfield, Kan \cite{Bousfield-Kan1}, Chap. {\bf I}, 3.2) given by $\Delta^k = \Delta[k]$.
The cosimplicial structure maps of $\Delta \tdot$ yield the simplicial structure maps of the simplicial set $hom(X \tdot,Y \tdot)$.

We will also consider the cosimplicial space denoted $sk_n \Delta \tdot$, given in dimension $k$ by $sk_n(\Delta[k])$, with obvious coface and codegeneracy maps.
Observe that $sk_n \Delta \tdot$ injects into $\Delta \tdot$ .

\begin{defn}
The functors $tot$, $tot^n$: $sSets^\Delta \lra sSets$ from cosimplicial spaces to spaces are defined by
\begin{center}
$tot X \tdot = hom(\Delta \tdot,X \tdot)$

$tot^n X \tdot = hom(sk_n \Delta \tdot ,X \tdot)$
\end{center}
\end{defn}

\bigskip

From this definition, it follows that there is a tower of maps 
\begin{center}
$tot X \tdot \ra ... \ra tot^{n+1} X \tdot \ra tot^n X \tdot \ra ... \ra tot^0 X \tdot $ $\cong$ $ X_0$
\end{center}

\bigskip

and $tot X \tdot $ $\cong$ $ {\lim}^{n \in {\mathbb N}} (tot^n X \tdot)$
where ${\mathbb N}$ = $\{0, 1, 2, ...\}$

\bigskip

\subsection{Skeleton and coskeleton}
\label{subsecCosk}
To investigate basic properties of $tot$ and $tot^n$, we turn our attention to the notion of skeleton and coskeleton.

Fix ${\cat C}$ a complete and cocomplete category. We investigate the skeleton and coskeleton functors for ${\cat C}^\Delta$.
Since we are not particularly interested in skeleton and coskeleton for ${\cat C}^{\Delta^{op}}$, we will skip a presentation of that.

\bigskip

Denote $\Delta_n$ the full subcategory of $\Delta$ with objects $\underline{0}, \underline{1}, ... \underline{n}$, and $i_n : \Delta_n \ra \Delta$ the embedding.

There are two pairs of adjoint functors

\begin{center}
\[ {\colim}^{i_n} : {\cat C}^{\Delta_n} \xymatrix{ \ar@<.5ex>[r] & \ar@<.5ex>[l]} {\cat C}^\Delta : {i_n^*} \]

\[ {i_n^*} : {\cat C}^\Delta \xymatrix{ \ar@<.5ex>[r] & \ar@<.5ex>[l]} {\cat C}^{\Delta_n} : {\lim}^{i_n} \]
\end{center}
where ${\colim}^{i_n}$, ${\lim}^{i_n}$ are the left, respectively the right Kan extensions along the functor $i_n$.

The n-th skeleton $sk_n X \tdot$ of a cosimplicial object $X \tdot \in Ob {\cat C}^\Delta$ is defined as 
\begin{center}
$sk_n X \tdot = {\colim}^{i_n} ({i_n^*} X \tdot)$.
\end{center}

The n-th coskeleton $cosk^n X \tdot$ of a cosimplicial object $X \tdot \in Ob {\cat C}^\Delta$ is defined as 
\begin{center}
$cosk^n X \tdot = {\lim}^{i_n} ({i_n^*} X \tdot)$.
\end{center}

\bigskip

We have the following useful characterization of the skeleton and the coskeleton:

\begin{prop}
\label{propSkCosk1}
If $X \tdot$ is an object in ${\cat C}^\Delta$, then:
\begin{list}{-}{}
 \item for $k \le n$, $(sk_n X \tdot)^k$ $\cong$ $X^k$ $\cong$ $(cosk^n X \tdot)^k$ under the adjunction maps $sk_n X \tdot \ra X \tdot \ra cosk^n X \tdot$
 \item denoting latching spaces by $L_k(-)$, the natural map \\
$L_{n + 1}(X\tdot) \ra (sk_n X\tdot)^{n + 1}$ is an isomorphism
 \item for $k \ge n+2$, the latching map $L_k(sk_n X\tdot) \ra (sk_n X\tdot)^k$ is an isomorphism
 \item denoting matching spaces by $M_k(-)$, the natural map \\
$(cosk^n X\tdot)^{n + 1} \ra M_{n + 1}(X\tdot)$ is an isomorphism
 \item for $k \ge n+2$, the matching map $(cosk^n X\tdot)^k \ra M_k(cosk^n X\tdot)$ is an isomorphism
\end{list}
These properties characterize uniquely $sk_n X \tdot$ and $cosk^n X \tdot$, for $X \tdot \in Ob {\cat C}^\Delta$.
\end{prop}

The proof for this proposition is standard and we omit it.

\bigskip
\bigskip

Observe that the functors $sk_n$, $cosk^n$ form an adjoint pair:

\begin{prop}
\label{propSkCosk2}
For $X \tdot$, $Y \tdot \in Ob {\cat C}^\Delta$, there is a natural isomorphism
\begin{center}
$Hom_{{\cat C}^\Delta} (sk_n X \tdot, Y \tdot)$ $\cong$ $Hom_{{\cat C}^\Delta} (X \tdot, cosk^n Y \tdot)$
\end{center}
\end{prop}

\begin{proof}
We have that
\begin{center}
$Hom_{{\cat C}^\Delta} ({\colim}^{i_n} ({i_n^*} X \tdot), Y \tdot)$ $\cong$ $Hom_{{\cat C}^{\Delta_n}} ({i_n^*} X \tdot, {i_n^*} Y \tdot) $

$\cong$ $Hom_{{\cat C}^\Delta} (X \tdot, {\lim}^{i_n} ({i_n^*} Y \tdot))$
\end{center}
\end{proof}

We now specialize to the case when ${\cat C}$ = $sSets$.

An easy computation using Proposition \ref{propSkCosk1} shows that the $n$-skeleton of the cosimplicial standard space $\Delta \tdot \in Ob(sSets^\Delta)$ is given in cosimplicial dimension $k$ by $sk_n(\Delta[k])$, with coface and codegeneracy maps induced from $\Delta \tdot$ .

Thus, the notation we used in Section \ref{subsectottotn} for $sk_n \Delta \tdot$ is consistent with the notion of skeleton in $sSets^\Delta$ that we introduced in the current Section.
As an aside, please note that it is not true in general that $(sk_n X \tdot)^k$ $\cong$ $sk_n (X^k)$ for a cosimplicial space $X \tdot$.

We conclude with the following useful Proposition:

\begin{prop}
\label{propTotCotot1}
The natural maps below are isomorphisms, for $X \tdot$ a cosimplicial space:
\begin{center}
$\xymatrix{
tot (cosk^n X \tdot) \ar[r]^\cong &
tot^n (cosk^n X \tdot) &
tot^n (X \tdot) \ar[l]_(.4)\cong
}$
\end{center}
\end{prop}

\begin{proof}
The functors $sk_n$, $cosk^n$ are adjoint (Proposition \ref{propSkCosk2}): for cosimplicial spaces $X \tdot$, $Y \tdot$ we have that
\begin{center}
$Hom_{sSets^\Delta} (sk_n X \tdot, Y \tdot)$ $\cong$ $Hom_{sSets^\Delta} (X \tdot, cosk^n Y \tdot)$
\end{center}
Observe that for a cosimplicial space $X \tdot$ and for a simplicial set $K$ we have that
\begin{center}
$sk_n (X \tdot \times K)$ $\cong$ $sk_n (X \tdot) \times K$
\end{center}
and from this we deduce that for cosimplicial spaces $X \tdot$, $Y \tdot$ we have a canonical isomorphism
\begin{equation}
\label{eqnhomAdjuncskcosk}
hom (sk_n X \tdot, Y \tdot) \cong hom (X \tdot, cosk^n Y \tdot)
\end{equation}

\bigskip

We are ready to prove that the three terms in Propostion \ref{propTotCotot1} are isomorphic.
From the definitions, the first term is just $hom(\Delta \tdot, cosk^n X \tdot)$, the second term is given by $hom(sk_n \Delta \tdot , cosk^n X \tdot) \cong hom(\Delta \tdot, cosk^n cosk^n X \tdot)$, and the third term is just $hom(sk_n \Delta \tdot , X \tdot)$.
They are all three isomorphic because of the isomorphism (\ref{eqnhomAdjuncskcosk}), and because $cosk^n X \tdot \cong cosk^n cosk^n X \tdot$ . 
\end{proof}

\bigskip

\subsection{$tot$ and $tot^n$ in terms of homotopy limits}
\label{subsecTotTotnHomotProp}

In this Subsection, we will investigate homotopical properties of the functors $tot$ and $tot^n$.

We will show that for a Reedy fibrant cosimplicial space $X \tdot$ we have natural weak equivalences

\begin{equation}
\label{eqnTotLim3}
tot X \tdot \simeq {\bf R}{\lim}^\Delta X \tdot 
\end{equation}

\begin{equation}
\label{eqnTotLim2}
tot^n X \tdot \simeq {\bf R}{\lim}^{\Delta_n} (i_n^* X\tdot)
\end{equation}
and the tower of maps
\begin{center}
$tot X \tdot \ra ... \ra tot^{n+1} X \tdot \ra tot^n X \tdot \ra ... \ra tot^0 X \tdot $ $\cong$ $ X^0$
\end{center}
is a tower of fibrant spaces and fibrations with
\begin{equation} 
\label{eqnTotLim4}
tot X \tdot \simeq {\bf R}{\lim}^{n \in {\mathbb N}} (tot^n X \tdot)
\end{equation}

\bigskip
\bigskip

We start by recalling the canonical isomorphism with $X \in Ob(sSets)$ and $Y \tdot \in Ob(sSets^\Delta)$
\begin{center}
$ Hom_{sSets^\Delta}(X \times \Delta \tdot, Y \tdot) = Hom_{sSets}(X , hom(\Delta \tdot, Y \tdot)) $
\end{center}
It follows that there is an adjoint pair of functors $F_\Delta$ = $(- \times \Delta \tdot)$, $tot$ = $hom(\Delta \tdot, -)$

\bigskip
\begin{center}
$
F_\Delta : \phantom{pp} sSets \xymatrix{ \ar@<.5ex>[r] & \ar@<.5ex>[l] } sSets^\Delta  \phantom{pp} : tot
$
\end{center}

\bigskip

The functor $F_\Delta$ = $(- \times \Delta \tdot)$ takes weak equivalences to weak equivalences.
Since $\Delta \tdot$ is Reedy cofibrant, it is easy to prove that $F_\Delta$ takes cofibrations to Reedy cofibrations.

It follows that the adjoint pair ($F_\Delta$, $tot$) is a Quillen adjoint pair, so the adjoint pair of left resp. right derived functors denoted
\bigskip
\begin{center}
${\bf L}F_\Delta$: ${\bf ho}(sSets) \xymatrix{ \ar@<.5ex>[r] & \ar@<.5ex>[l] } {\bf ho}(sSets^\Delta )$  : ${\bf R}tot$
\end{center}
exists.

Since $\Delta[n]$ is contractible for all $n$, it follows that ${\bf L}F_\Delta$ = ${\bf L}c$, where $c$ denotes the constant functor $c$: $sSets \lra sSets^\Delta$.
Since ${\bf R}tot$ is right adjoint to ${\bf L}c$, it follows that ${\bf R}tot$ = ${\bf R}\lim^\Delta$.

\bigskip
\bigskip

Because of the Quillen adjunction ($F_\Delta, tot$), it follows that for $X \tdot$ Reedy fibrant $tot X \tdot$  computes ${\bf R}tot (X \tdot)$:
\begin{center}
$tot X \tdot$ $\simeq$ ${\bf R}{\lim}^\Delta X \tdot$
\end{center}
and this proves Equation (\ref{eqnTotLim3}).

\bigskip
\bigskip

It is easy to see using Proposition \ref{propSkCosk1} that $cosk^n X\tdot$ is also Reedy fibrant.
By Proposition \ref{propTotCotot1}, we also have
\begin{center}
$tot^n X \tdot \simeq {\bf R}{\lim}^\Delta (cosk^n X\tdot)$
\end{center}

The adjoint pair below is a Quillen adjoint pair:
\begin{center}
$ {i_n^*} : sSets^\Delta \xymatrix{ \ar@<.5ex>[r] & \ar@<.5ex>[l]} sSets^{\Delta_n} : {\lim}^{i_n} $
\end{center}

Therefore ${\lim}^{i_n} A \tdot$ $\simeq$ ${\bf R}{\lim}^{i_n} A$, for $A \tdot$ a Reedy fibrant object of $sSets^{\Delta_n}$.
Since $X \tdot$ is Reedy fibrant, $i_n^* X\tdot$ is Reedy fibrant in $sSets^{\Delta_n}$.
It follows that
\begin{center}
$cosk^n X\tdot$ = ${\lim}^{i_n} (i_n^* X\tdot)$ $\simeq$ $({\bf R}{\lim}^{i_n}) (i_n^* X\tdot)$
\end{center}

\bigskip
Dwyer, Kan and Hirschhorn \cite{Dwyer-Kan-Hirschhorn} 58.3 prove a composability property for relative homotopy limits.
From their general result we deduce that the natural map

\begin{center}
$\xymatrix{
{\bf R}{\lim}^\Delta \comp {\bf R}{\lim}^{i_n} (i_n^* X\tdot) \ar[r]^(.57)\simeq &
{\bf R}{\lim}^{\Delta_n} (i_n^* X\tdot)
}$
\end{center}
is a weak equivalence.

This leads to the natural weak equivalence

\begin{center}
$tot^n X \tdot$ $\simeq$ ${\bf R}{\lim}^{\Delta_n} (i_n^* X\tdot)$
\end{center}
and this proves Equation (\ref{eqnTotLim2}).

\bigskip
\bigskip

One can easily show that the natural weak equivalences (\ref{eqnTotLim3}), (\ref{eqnTotLim2}) are compatible under the natural maps 
\begin{center}
$tot X \tdot \ra tot^n X \tdot$ and ${\bf R}{\lim}^\Delta X \tdot \ra {\bf R}{\lim}^{\Delta_n} (i_n^* X\tdot)$  

$tot^{n+1} X \tdot \ra tot^n X \tdot$ and ${\bf R}{\lim}^{\Delta_{n+1}} (i_{n+1}^* X \tdot) \ra {\bf R}{\lim}^{\Delta_n} (i_n^* X\tdot)$  
\end{center}

\bigskip
\bigskip

Again by Proposition \ref{propSkCosk1}, if $X \tdot$ is Reedy fibrant, then
\begin{center}
$cosk^{n + 1} X\tdot \xymatrix{\ar@{->>}[r] &} cosk^n X\tdot$
\end{center}
is a Reedy fibration of Reedy fibrant objects.
Furthermore, by the Quillen adjunction ($F_\Delta, tot$), $tot$ sends Reedy fibrations to fibrations, so we get that
\begin{center}
$tot^{n + 1} X\tdot \xymatrix{\ar@{->>}[r] &} tot^n X\tdot$
\end{center}
is a fibration of fibrant simplicial sets.
Putting all this together, if $X \tdot$ is Reedy fibrant, the tower of maps
\begin{center}
$tot X \tdot \ra ... \ra tot^{n+1} X \tdot \ra tot^n X \tdot \ra ... \ra tot^0 X \tdot $ $\cong$ $ X^0$
\end{center}
is a tower of fibrant spaces and fibrations, and 
\begin{center}
$tot X \tdot $ $\cong$ $ {\lim}^{n \in {\mathbb N}} (tot^n X \tdot)$ $\simeq$ $ {\bf R}{\lim}^{n \in {\mathbb N}} (tot^n X \tdot)$
\end{center}
which proves Equation (\ref{eqnTotLim4}).

\bigskip
\bigskip

\subsection{Pointed $tot$ and $tot_n$}
\label{subsecpointedtotcotot}
In this Section, we investigate the functors $tot_*$, $tot^n_*$ : $sSets_*^\Delta \lra sSets_*$ , which are pointed versions of $tot$ and $tot^n$.

For $X \tdot \in Ob(sSets^\Delta)$ and $Y \tdot \in Ob(sSets_*^\Delta)$, we define the pointed simplicial hom-set $hom_*(X \tdot, Y \tdot) \in Ob(sSets_*)$ by choosing basepoints in the following way:

\begin{center}
$hom_*(X \tdot,Y \tdot)_k$ = cosimplicial maps $X \tdot \times \Delta [k] \ra Y \tdot$

$\phantom{foooooooooooooo}$ with basepoint $\phantom{fo} X \tdot \times \Delta [k] \ra *$ 
\end{center}

\bigskip

The functors $tot_*$, $tot^n_*$: $sSets_*^\Delta \lra sSets_*$ from cosimplicial pointed spaces to pointed spaces are defined by
\begin{center}
$tot_* X \tdot = hom_*(\Delta \tdot,X \tdot)$

$tot^n_* X \tdot = hom_*(sk_n \Delta \tdot ,X \tdot)$
\end{center}

\bigskip

Keeping track of basepoints, one easily sees that for $X \tdot \in Ob(sSets_*^\Delta)$ the natural maps below are isomorphisms:
\begin{center}
$\xymatrix{
tot_* (cosk^n X \tdot) \ar[r]^\cong &
tot^n_* (cosk^n X \tdot) &
tot^n_* (X \tdot) \ar[l]_(.4)\cong
}$
\end{center}

\bigskip

We have a left adjoint functor to $tot_*$:

\begin{center}
\bigskip
$F_{\Delta *} : \phantom{pp} sSets_* \xymatrix{ \ar@<.5ex>[r] & \ar@<.5ex>[l] } sSets_*^\Delta  \phantom{pp} : tot_*$
\bigskip
\end{center}
given by $F_{\Delta *}(X)$ = $X \times \Delta \tdot / * \times \Delta \tdot $.

The functor $tot_*$ takes Reedy fibrations to fibrations and equivalences between Reedy fibrant objects to equivalences, because $tot$ does so.
It follows that the pair ($F_{\Delta *}$, $tot_*$) is a Quillen adjoint pair, so the adjoint pair of left resp. right derived functors exists:
\bigskip
\begin{center}
${\bf L}F_{\Delta *}$: ${\bf ho}(sSets_*) \xymatrix{ \ar@<.5ex>[r] & \ar@<.5ex>[l] } {\bf ho}(sSets_*^\Delta )$  : ${\bf R}tot_*$
\end{center}

The functor $F_{\Delta *}$ takes weak equivalences to weak equivalences, and ${\bf L}F_{\Delta *}$ = ${\bf L}c$, where $c$ : $sSets_* \lra sSets_*^\Delta$ is the constant functor.
From this observation, we conclude that ${\bf R}tot_*$ = ${\bf R}\lim^\Delta$. 

\bigskip

The proofs of Section \ref{subsecTotTotnHomotProp} yield easily that for a Reedy fibrant cosimplicial pointed space $X \tdot \in Ob(sSets_*^\Delta)$ we have natural pointed weak equivalences

\begin{equation}
\label{eqnTotLim3*}
tot_* X \tdot \simeq {\bf R}{\lim}^\Delta X \tdot 
\end{equation}

\begin{equation}
\label{eqnTotLim2*}
tot^n_* X \tdot \simeq {\bf R}{\lim}^{\Delta_n} (i_n^* X\tdot)
\end{equation}
and the tower of maps
\begin{center}
$tot_* X \tdot \ra ... \ra tot^{n+1}_* X \tdot \ra tot^n_* X \tdot \ra ... \ra tot^0_* X \tdot $ $\cong$ $ X^0$
\end{center}
is a tower of fibrant pointed spaces and fibrations with
\begin{equation} 
\label{eqnTotLim4*}
tot_* X \tdot \simeq {\bf R}{\lim}^{n \in {\mathbb N}} (tot^n_* X \tdot)
\end{equation}

\bigskip
\bigskip

\subsection{tot and cotot}
\label{subsecTotCotot}

In this Subsection, let ${\cat M}$ be a model category. 
We extend the results in the previous subsections to the case of the model category ${\cat M}$. 

For $X \tdot$ in ${\bf ho}({\cat M}^\Delta)$, we define ${\bf tot} X \tdot$ and ${\bf tot}^n X \tdot$ as objects in ${\bf ho}{\cat M}$:

\begin{center}
${\bf tot} X \tdot = {\bf R}{\lim}^{\Delta} X \tdot$

${\bf tot}^n X \tdot = {\bf R}{\lim}^{\Delta_n} {i_n^*}X \tdot$
\end{center}

\bigskip

Note the essential difference between $tot$, $tot^n$ and ${\bf tot}$, ${\bf tot}^n$: the first ones live in the model category $sSets$, whereas the last ones live in the homotopy category ${\bf ho}{\cat M}$ of a model category.

\bigskip

In the case ${\cat M}$ = $sSets$, if $X \tdot$ is a Reedy fibrant cosimplicial space, then by (\ref{eqnTotLim3}),  (\ref{eqnTotLim2}) ${\bf tot} X \tdot \simeq tot X \tdot$, and ${\bf tot}^n X \tdot \simeq tot^n X \tdot$.

\bigskip

Dually, for $Y .$ in ${\bf ho}({\cat M}^{\Delta^{op}})$ we define ${\bf cotot} Y.$ and ${\bf cotot}_n Y.$ as objects in ${\bf ho}{\cat M}$:

\begin{center}
${\bf cotot} Y. = {\bf L}{\colim}^{\Delta^{op}} Y.$

${\bf cotot}_n Y. = {\bf L}{\colim}^{\Delta_n^{op}} {(i_n^{op})^* Y.}$
\end{center}

\bigskip

We get natural maps 
\begin{center}
${\bf tot} X \tdot \ra ... \ra {\bf tot}^{n+1} X \tdot \ra {\bf tot}^n X \tdot \ra ... \ra {\bf tot}^0 X \tdot $ $\cong$ $ X^0$
\end{center}
We can see this sequence of maps as an object in ${\bf ho}({\cat M}^{\overline {\mathbb N}})$, where ${\overline {\mathbb N}}$ has objects $0$, $1$, $2$, ... , $\infty$ and unique maps $i \ra j$ for $i \ge j$.

Dually we get natural maps
\begin{center}
${\bf cotot} Y. \la ... \la {\bf cotot}_{n+1} Y. \la {\bf cotot}_n Y. \la ... \la {\bf cotot}_0 Y. $ $\cong$ $ Y_0$
\end{center}
and we can see this sequence of maps as an object in ${\bf ho}({\cat M}^{{\overline {\mathbb N}}^{op}})$.

For convenience, we denote ${\bf tot}^\infty = {\bf tot}$, ${\bf cotot}_\infty = {\bf cotot}$.

Consequences of Theorem \ref{thmC2} and (\ref{eqnTotLim3}), (\ref{eqnTotLim2}), (\ref{eqnTotLim4}) are the two theorems below:

\begin{thm}
\label{thmTotCototFunCom}
For all $n \in {\overline {\mathbb N}}$, there are natural weak equivalences, compatible with maps $n_1 \ra n_2$ in ${\overline {\mathbb N}}$
\begin{center}
${\bf Hom}(W, {\bf tot}^n X \tdot) \cong ({\bf R}tot^n) {\bf Hom}(W, X \tdot)$

${\bf Hom}({\bf cotot}_n Y., W) \cong ({\bf R}tot^n) {\bf Hom}(Y., W)$  
\end{center}
\end{thm}

\begin{thm}
The natural maps ${\bf tot} \ra {\bf tot}^n$ , ${\bf cotot}_n \ra {\bf cotot} $ yield isomorphisms
\begin{center}
${\bf tot} X \tdot $ $\cong$ $ {\bf R}{\lim}^{n \in {\mathbb N}} {\bf tot}^n X \tdot$

${\bf cotot} Y. $ $\cong$ $ {\bf L}{\colim}^{n \in {{\mathbb N}}^{op}} {\bf cotot}_n Y.$
\end{center}
\end{thm}

\bigskip

The upshot of using $tot_*$ and $tot_*^n$ is the following pointed version of Thm \ref{thmTotCototFunCom}:

\begin{thm}
\label{thmTotCommPointed}
For a {\bf pointed} model category ${\cat M}$, for all $n \in {\overline {\mathbb N}}$, there are natural weak equivalences compatible with maps $n_1 \ra n_2$ in ${\overline {\mathbb N}}$
\begin{center}
${\bf Hom}_*(W, {\bf tot}^n X \tdot) \cong ({\bf R}tot^n_*) {\bf Hom}_*(W, X \tdot)$

${\bf Hom}_*({\bf cotot}_n Y., W) \cong ({\bf R}tot^n_*) {\bf Hom}_*(Y., W)$  
\end{center}
\end{thm}

\bigskip
\bigskip
\bigskip
\section{The Spectral Sequence}
\label{secSpSeq}

\bigskip
\bigskip

In this Section we present the homotopy spectral sequence of a cosimplicial object in a pointed model category.
In the second part, we define the dual (homology) spectral sequence.

\subsection{The homotopy spectral sequence of a cosimplicial object}
\label{subsecHomotopySpSeq}

Let ${\cat M}$ be a pointed model category.

Fix any representative $\mathcal{H}om_*($-$,$-$)$ in $sSets_*^{{\cat M}^{op}\times{\cat M}}$ of the pointed function complex ${\bf Hom}_*($-$,$-$)$.

For $X \tdot$ a cosimplicial object in ${\cat M}$ and $W$ an object in ${\cat M}$, form the cosimplicial pointed space 
\begin{center}
$ \mathcal{H}om_*(W, X \tdot) $

\bigskip
\end{center}
Since $\mathcal{H}om_*(W, X \tdot)$ is not necessarily Reedy fibrant in $sSets_*^\Delta$, we choose a canonical Reedy fibrant replacement in $sSets^\Delta$ denoted $\mathcal{H}om_*(W, X \tdot)'$.

\bigskip

Define the homotopy spectral sequence as
\begin{equation}
E_r^{s,t}(W, X \tdot) = E_r^{s,t}(\mathcal{H}om_*(W, X \tdot)') \hspace{.5in} 
(t \ge s \ge 0)
\end{equation}
where $E_r^{s,t}(\mathcal{H}om_*(W, X \tdot)')$ denotes the homotopy spectral sequence of the cosimplicial pointed space $\mathcal{H}om_*(W, X \tdot)'$ (\cite{Bousfield-Kan1}, Chap. {\bf X} Sec. 6).

This spectral sequence computes the homotopy groups of 
\begin{center}
$tot_* \mathcal{H}om_*(W, X \tdot)' $ $\simeq$ $ {\bf Hom}_*(W, {\bf tot}X \tdot)$ 
\end{center}
by Theorem \ref{thmTotCommPointed}.

\bigskip

There is actually a better way to understand this spectral sequence, using the ${\bf tot}$-tower of $X \tdot$.
The homotopy spectral sequence $E_r^{s,t}(W, X \tdot)$ is just the spectral sequence of the $tot$-tower of pointed spaces
\begin{center}
... $\ra tot_*^{n+1}\mathcal{H}om_*(W, X \tdot)' \ra tot_*^n \mathcal{H}om_*(W, X \tdot)' \ra$ ... $\ra tot_*^0 \mathcal{H}om_*(W, X \tdot)'$ 
\end{center}

By Theorem \ref{thmTotCommPointed}, this $tot$-tower is weakly equivalent in $sSets^{\bf N}_*$ to
\begin{center}
... $\ra {\bf Hom}_*(W, {\bf tot}^{n+1} X \tdot) \ra {\bf Hom}_*(W, {\bf tot}^n X \tdot) \ra$ ... $\ra {\bf Hom}_*(W, {\bf tot}^0 X \tdot)$ 
\end{center}
so we can say that $E_r^{s,t}(W, X \tdot)$ is the spectral sequence of the (relative homotopy groups with coefficients in $W$ of the) ${\bf tot}$-tower of $X \tdot$.

\bigskip

\subsection{The homology spectral sequence of a simplicial object}

For $Y.$ a simplicial object in ${\cat M}$ and $W$ an object in ${\cat M}$, form the cosimplicial pointed space

\bigskip
\begin{center}
$ \mathcal{H}om_*(Y., W) $
\end{center}
\bigskip

Choose a canonical Reedy fibrant replacement in $sSets^\Delta$ for $\mathcal{H}om_*(Y., W)$, denoted $\mathcal{H}om_*(Y., W)'$

\bigskip

Define the homology spectral sequence as
\begin{equation}
E_r^{s,t}(Y., W) = E_r^{s,t}(\mathcal{H}om_*(Y., W)') \hspace{.5in} 
(t \ge s \ge 0)
\end{equation}
where $E_r^{s,t}(\mathcal{H}om_*(Y, W \tdot)')$ denotes the homotopy spectral sequence of the cosimplicial pointed space $\mathcal{H}om_*(Y, W \tdot)'$.

This spectral sequence computes the homotopy groups of ${\bf Hom}_*({\bf cotot Y.}, W)$
and is the same as the homotopy spectral sequence of the tower in ${\bf ho}(sSets^{\bf N}_*)$
\begin{center}
... $\ra {\bf Hom}_*({\bf cotot}_{n+1} Y., W) \ra {\bf Hom}_*({\bf cotot}_n Y., W) \ra$ ... $\ra {\bf Hom}_*({\bf cotot}_0 Y., W)$ 
\end{center}

\bigskip

\bigskip
                % Local Variables: 
                % mode: latex 
                % tex-command: amslatex-command
                % tex-dvi-view-command: tex-oneside-view-command
                % tex-start-of-header: "\\special"
                % tex-end-of-header: "\\maketitle"
                % abbrev-mode: t
                % End:

%% file: appa.tex
\chapter{Homotopy Limits and Simplicial Objects}
\label{app1}

\bigskip
\bigskip

In this Appendix, for a model category, we prove that two homotopic maps of cosimplicial objects coincide after taking the homotopy limit.

\bigskip

Here is what we do in each of the two Sections:

In Section \ref{secHomCosMaps}, we recall the homotopy relation between  maps of (co)simplicial objects in a category.
The most illustrative example of homotopic maps is: if $f$, $g$ are two triple maps from a triple ($R, \nu, M$) to a triple ($T, \nu, M$), then the maps of the associated cosimplicial objects
\begin{center}
${\cat R} \tdot (X) \xymatrix{ \ar@<.5ex>[r]^f \ar@<-.5ex>[r]_g &} {\cat T} \, \tdot (X)$
\end{center}
are naturally homotopic (Prop. \ref{propMapsOfTriples}).

In Section \ref{secHomCosMaps2}, we specialize to the case of (co)simplicial objects in a model category.
We prove the main result of this Appendix (Thm. \ref{thmHomot2}), namely that two homotopic maps of cosimplicial objects 
\begin{center}
$f$ $\simeq$ $g$ : $X \tdot \lra Y \tdot$ 
\end{center}
have the same effect on the homotopy limit of these cosimplicial objects (regard the cosimplicial objects as diagrams indexed by the category $\Delta$):

\bigskip
\begin{center}
${\mathbf R}\lim^\Delta f$ = ${\mathbf R}\lim^\Delta g$
\end{center}

\bigskip
\bigskip
\bigskip

\section{The Homotopy Relation}
\label{secHomCosMaps}

\bigskip
\bigskip

Let ${\cat C}$ be a complete and cocomplete category.
\subsection{$X. \tens K.$ and ${(Y \tdot)}^{K.}$}
Recall that the following two functors define an action of simplicial sets on ${\cat C}^{\Delta^{op}}$, and a coaction of simplicial sets on ${\cat C}^\Delta$:

\bigskip

\begin{center}
${\cat C}^{\Delta^{op}} \times sSets \lra {\cat C}^{\Delta^{op}}$

($X.$ , $K.$) $\mapsto$ $X. \tens K.$ 

\end{center}
where $(X. \tens K.)_n = \sqcup_{K_n} X_n$, and dually

\begin{center}
${\cat C}^\Delta \times sSets^{op} \lra {\cat C}^\Delta$

($Y \tdot$ , $K.$) $\mapsto$ $(Y \tdot)^{K.}$

\end{center}
where $({(Y \tdot)}^{K.})^n = \times_{K_n} Y^n$.

\bigskip
\bigskip

\subsection{The homotopy relation}
If $f$, $g$ : $X. \xymatrix{ \ar@<.5ex>[r] \ar@<-.5ex>[r] &} Y.$ are two maps in ${\cat C}^{\Delta^{op}}$, we say that $f$ is homotopic to $g$ (write $f \simeq g$ ) if there exists a map $H$ that makes the diagram below commute:

\begin{center}
$\xymatrix{ X. \cong X. \tens \Delta[0] \ar[d]_{X. \tens d^0} \ar[dr]^f & \\ 
	X. \tens \Delta[1] \ar[r]^H &
	Y. \\
	X. \cong X. \tens \Delta[0] \ar[u]^{X. \tens d^1} \ar[ur]_g &
}$  
\end{center}

\bigskip

Dually, if $f$, $g$ : $X \tdot \xymatrix{ \ar@<.5ex>[r] \ar@<-.5ex>[r] &} Y \tdot$ are two maps in ${\cat C}^\Delta$, then $f \simeq g$ if there exists a map $H$ that makes the following diagram commute:
\begin{center}
$\xymatrix{ & Y \tdot \cong (Y \tdot)^{\Delta[0]} \\ 
	X \tdot \ar[r]^H \ar[ur]^f \ar[dr]_g &
	(Y \tdot)^{\Delta[1]} \ar[u]_{Y^{d^0}} \ar[d]^{Y^{d^1}}  \\
	& Y \tdot \cong (Y \tdot)^{\Delta[0]} 
}$  
\end{center}

\bigskip
\bigskip

Let $F$ : ${\cat C}_1 \ra {\cat C}_2$ be a functor.
The following two propositions are straightforward:

\begin{prop}
\label{propBKnew1}
If $f$ $\simeq$ $g$ : $X. \lra Y.$ are two homotopic maps in ${\cat C}_1^{\Delta^{op}}$, then $F(f)$ $\simeq$ $F(g)$ in ${\cat C}_2^{\Delta^{op}}$
\end{prop}
and dually
\begin{prop}
\label{propBKnew2}
If $f$ $\simeq$ $g$ : $X \tdot \lra Y \tdot$ are two homotopic maps in ${\cat C}_1^\Delta$, then $F(f)$ $\simeq$ $F(g)$ in ${\cat C}_2^\Delta$.
\end{prop}

\bigskip

Let us fix some more terminology, before going to the next subsection.
We say that a (co)simplicial map $f$ is a homotopy equivalence if there's a (co)simplicial map $g$ in the other direction such that $gf \simeq id$, $fg \simeq id$.
If $X. \ra X_{-1}$ is an augmented simplicial object in ${\cat C}$, and if the map of simplicial objects 

\begin{center}
$X. \ra cX_{-1}$
\end{center}
to the constant simplicial object is a homotopy equivalence, then we say that the whole augmented simplicial object ($X. \ra X_{-1}$) is contractible.
If $Y^{-1} \ra Y \tdot$ is an augmented cosimplicial object, and if 

\begin{center}
$cY^{-1} \ra Y \tdot$ 
\end{center}
is a homotopy equivalence, we say that the whole augmented cosimplicial object $(Y^{-1} \ra Y \tdot)$ is contractible.

\bigskip
\bigskip

\subsection{Examples arising from triples and cotriples}
Several typical examples of homotopic maps come from contructions that involve triples and cotriples.

For any triple ($R, \nu, M$) on ${\cat C}$
\begin{center}
$ \xymatrix {
	X \ar[r]^\nu &
	RX &
	R^2X \ar[l]_M
	} $
\end{center}
one forms the natural augmented cosimplicial resolution ${\cat R} \tdot (X)$ of $X$

\begin{equation}
\label{eqnBC5}
\overline{{\cat R}} \; \tdot (X) : \phantom{foo} X \lra {\cat R} \tdot (X)
\end{equation}
given by
\begin{center}
\[ \xymatrix{ 
	{\overline{{\cat R}} \; \tdot} (X): &  
	X \ar[r] & 
	RX \ar[r] \ar@<1ex>[r] & 
	R^2X \ar@<1ex>[l] \ar[r] \ar@<1ex>[r] \ar@<2ex>[r] & 
	R^3X \;\; ... \ar@<1ex>[l] \ar@<2ex>[l]
	}
\]
\end{center}

\bigskip

It is defined by the following formulas:

\begin{center}
${\overline {\cat R} \,}^n (X) = R^{n+1}X \phantom{poo} (n \ge -1)$

${\cat R}^n (X) = {\overline {\cat R} \,}^n (X) \phantom{poo} (n \ge 0)$

$d^i : {\overline {\cat R} \,}^{n-1}(X) \ra {\overline {\cat R} \,}^n(X), \phantom{po} d^i = R^i \nu R^{n-i} \phantom{poo} (n \ge 0, n \ge i \ge 0)$

$s^i : {\overline {\cat R} \,}^{n+1}(X) \ra {\overline {\cat R} \,}^n(X) , \phantom{po} s^i = R^i M R^{n-i} \phantom{poo} (n \ge 0, n \ge i \ge 0)$
\end{center}

\bigskip

Observe that by (\ref{eqnBC5}) we can regard $\overline{{\cat R}} \; \tdot (X)$ as a natural map from the constant cosimplicial object $cX$ to the cosimplicial object ${\cat R} \tdot (X)$.
The following Proposition yields an important example of contractible cosimplicial resolutions:

\begin{prop}
\label{propContractCot}
The augmented cosimplicial objects $\overline{{\cat R}} \; \tdot (RX)$ and $R\overline{{\cat R}} \; \tdot (X)$ are naturally contractible.
\end{prop}

We don't include a proof for this Proposition.

Instead, we choose to present in detail another important example of homotopic maps, that arises from considering maps of triples (refining an idea of Meyer, \cite{Meyer}).

To start, let us recall the definition of a map of triples.
Suppose we have two triples ($R, \nu, M$) and ($T, \nu, M$) on ${\cat C}$
\begin{center}
$ \xymatrix {
	X \ar[r]^\nu &
	RX &
	R^2X \ar[l]_M
	} $

$ \xymatrix {
	X \ar[r]^\nu &
	TX &
	T^2X \ar[l]_M
	} $

\bigskip
\end{center}
For simplicity, we use the same letters to denote the structure maps $\nu$, $M$ of the two triples.

A triple map $f$ from ($R, \nu, M$) to ($T, \nu, M$) by definition makes the two diagrams below commutative
\begin{center}
\bigskip
$ \xymatrix {
	X \ar[r]^\nu \ar[dr]_\nu &
	RX \ar[d]^f \\
	& TX
	} $

\bigskip

$ \xymatrix {
	RX \ar[dd]_f &
	R^2X \ar[l]_M \ar[d]^{Rf} \\
	& RTX \ar[d]^{fT} \\
	TX &
	T^2X \ar[l]^M
	} $
\end{center}

If $g$ is a second triple map from ($R, \nu, M$) to ($T, \nu, M$) then $fT \comp Rg = Tg \comp fR$, and we denote these compositions by $fg : R^2X \ra T^2X$.

Let's extend this notation: if $f_1$, ... , $f_n$ are triple maps from ($R, \nu, M$) to ($T, \nu, M$), denote by the juxtaposition $f_1f_2...f_n : R^nX \ra T^nX$ the composition 
\begin{center}
$f_1R^{n-1} \comp Tf_2R^{n-2} \comp ... T^{n-1}f_n$
\end{center}

\bigskip

We are ready for
\begin{prop}
\label{propMapsOfTriples}
If $f$, $g$ are two triple maps from ($R, \nu, M$) to ($T, \nu, M$), then the maps of cosimplicial objects
\begin{center}
${\cat R} \tdot (X) \xymatrix{ \ar@<.5ex>[r]^f \ar@<-.5ex>[r]_g &} {\cat T} \; \tdot (X)$
\end{center}
are naturally homotopic.
\end{prop}

\begin{proof}
We construct a homotopy
\begin{center}
$H$ : ${\cat R} \tdot (X) \lra ({\cat T} \; \tdot (X))^{\Delta[1]}$
\end{center}
by constructing $H^n$
\begin{center}
$H^n$ : $R^{n+1}(X) \lra \times_{\Delta[1]_n} T^{n+1} (X)$
\end{center}

Since $\Delta[1]_n = Hom_{sSets}(\Delta[n], \Delta[1]) = Hom_{\Delta}(\underline{n}, \underline{1})$, we need for any map $\underline{n} \ra \underline{1}$ to construct a map $R^{n+1}X \ra T^{n+1}X$.

If the map $\underline{n} \ra \underline{1}$ carries $0$, ..., $i$ to $0$ and carries $i+1$, ... $n$ to $1$ (where $-1 \leq i \leq n$), choose 
\begin{center}
$ff..gg$ : $R^{n+1}X \ra T^{n+1}X$
\end{center}
where in $ff..gg$ we have ($i+1$) consecutive copies of $f$ juxtaposed with ($n-i$) consecutive copies of $g$. 
\end{proof}

\bigskip
\bigskip
\bigskip

If ($U, \eps, \Delta$) is a cotriple on ${\cat C}$, construct the simplicial augmented object
\begin{center}
${\overline {\cat U}}.(X): \phantom{foo}
\xymatrix{
	{\cat U}.(X) \ar[r] &
	X
}
$
\end{center}
given by
\begin{center}
$\xymatrix{
	{\overline {\cat U} \,}^T \!\!\! . (X): & 
	... \;\; U^3X\ar[r] \ar@<1ex>[r] \ar@<2ex>[r] &
	U^2X \ar[r] \ar@<1ex>[r] \ar@<1ex>[l] \ar@<2ex>[l] &
	UX \ar[r] \ar@<1ex>[l] &
	X
}
$
\end{center}

The dual propositions are:
\begin{prop}
The augmented simplicial objects $\overline{{\cat U}}. (UX)$ and $U\overline{{\cat U}}. (X)$ are naturally contractible.
\end{prop}

\begin{prop}
If $f$, $g$ are two cotriple maps from the cotriple ($U, \eps, \Delta$) to the cotriple ($V, \eps, \Delta$), then the maps of simplicial objects
\begin{center}
${\cat U}. (X) \xymatrix{ \ar@<.5ex>[r]^f \ar@<-.5ex>[r]_g &} {\cat V}.(X)$
\end{center}
are naturally homotopic.
\end{prop}

\bigskip
\bigskip
\bigskip

\section{Effect of Homotopic Maps on Homotopy Limits}
\label{secHomCosMaps2}

\bigskip
\bigskip

In this Section we will prove:

\begin{thm}
\label{thmHomot2}
Suppose we have in ${\cat M}^\Delta$ two homotopic maps 
\begin{center}
$f$ $\simeq$ $g$ : $X \tdot \lra Y \tdot$ 
\end{center}
Then ${\bf R}\lim^\Delta f$ = ${\bf R}\lim^\Delta g$.
\end{thm}

Dually

\begin{thm}
\label{thmHomot1}
Suppose we have in ${\cat M}^{\Delta^{op}}$ two homotopic maps 
\begin{center}
$f$ $\simeq$ $g$: $X. \lra Y.$ 
\end{center}
Then ${\bf L}\colim^{\Delta^{op}} f$ = ${\bf L}\colim^{\Delta^{op}} g$.
\end{thm}

\begin{proof}
To prove Theorem \ref{thmHomot2} in general, because function complexes commute with homotopy limits (Chap. \ref{chapHSS}, Thm. \ref{thmC2}), it will be enough to prove Theorem \ref{thmHomot2} for the particular case when ${\cat M}$ = $sSets$.

This reduces to showing that if $Y \tdot$ is a cosimplicial space, then the map
\begin{center}
$\xymatrix{
	(Y \tdot)^{\Delta[0]} \ar[rr]^{(Y \tdot)^ {s^0}} & &
	(Y \tdot)^{\Delta[1]}
}$
\end{center}
induces an isomorphism
\begin{equation}
\label{eqnWhatToProveTot}
\xymatrix{
	{\bf R}\lim^\Delta (Y \tdot)^{\Delta[0]} \ar[rr]_\cong^{{\bf R}\lim^\Delta (Y \tdot)^ {s^0}} & &
	{\bf R}\lim^\Delta (Y \tdot)^{\Delta[1]}
}
\end{equation}

For a simplicial set $K.$ and a cosimplicial space $Y \tdot$ we could define a bicosimplicial space $[Y^K]\tdotdot$ by 
\begin{center}
$[Y^K]^{p, q}$ = $\times_{K_q} Y^p$
\end{center}
We have that $(Y \tdot)^{K.}$ = $diag [Y^K]\tdotdot$.

The reduction to the diagonal argument (Appendix \ref{app2}, Lemma \ref{lemmaBK1}) for the map
\begin{center}
$\xymatrix{
	[Y^{\Delta[0]}]\tdotdot \ar[rr]^{(Y \tdot)^ {s^0}} & &
	[Y^{\Delta[1]}]\tdotdot
}$
\end{center}
combined with the factorization ${\bf R}\lim^{(p,q) \in \Delta \times \Delta}$ = ${\bf R}\lim^{p \in \Delta}{\bf R}\lim^{q \in \Delta}$ quickly show that (\ref{eqnWhatToProveTot}) is true for cosimplicial spaces {\it if it is true for constant cosimplicial spaces}.

So, the proof reduces to showing that for a simplicial set $A$ (seen as a constant cosimplicial space) we have an isomorphism in ${\bf ho} sSets$
\begin{equation}
\label{eqnWhatToProveTot.1}
\xymatrix{
	{\bf R}\lim^\Delta (A^{\Delta[0]}) \ar[rr]_\cong^{{\bf R}\lim^\Delta (A^{s^0}) } & &
	{\bf R}\lim^\Delta (A^{\Delta[1]})
}
\end{equation}

We may assume that $A$ is a fibrant simplicial set.
The result now follows, because for a simplicial set $K$ and for a fibrant simplicial set $A$ we have that

\begin{center}
${\bf R}\lim^\Delta (A^K) \simeq hom(K, A)$
\end{center}
therefore the spaces in the equation (\ref{eqnWhatToProveTot.1}) are both weakly equivalent to $A$.
\end{proof}

\bigskip
\bigskip
\bigskip

                % Local Variables: 
                % mode: latex 
                % tex-command: amslatex-command
                % tex-dvi-view-command: tex-oneside-view-command
                % tex-start-of-header: "\\special"
                % tex-end-of-header: "\\maketitle"
                % abbrev-mode: t
                % End:

%% file: appb.tex
\chapter{The Reduction to the Diagonal Argument}
\label{app2}

\bigskip
\bigskip

In this Appendix we prove that the homotopy limit of a bicosimplicial object can be computed as the homotopy limit of the diagonal.
Dually, the homotopy colimit of a bisimplicial object is the homotopy colimit of the diagonal.

\bigskip
\bigskip

\section{Reduction to the Diagonal}

\bigskip
\bigskip

Here is what we want to prove:

\begin{lemma}[Reduction to the diagonal]
\label{lemmaBK1}
Let ${\cat M}$ be a model category.
If $X \tdotdot \in Ob({\cat M}^{\Delta \times \Delta})$ is a bicosimplicial object, then the natural map below is an ismorphisms in ${\bf ho} {\cat M}$:
\begin{center}
${\bf R}\lim^{\Delta \times \Delta} X \tdotdot$
$\xymatrix{\ar[r]^\cong &}$ 
${\bf R}\lim^\Delta \phantom{p} diag(X \tdotdot)$
\end{center}

Dually, if $X..$ is a bi-simplicial object, then
\begin{center}
${\bf L}\colim^{\Delta^{op}} \phantom{p} diag(X..)$
$\xymatrix{\ar[r]^\cong &}$ 
${\bf L}\colim^{\Delta^{op} \times \Delta^{op}} X..$ 
\end{center}
\end{lemma}

\bigskip

\begin{proof}
Based on the Cofinality Theorem \cite{Dwyer-Kan-Hirschhorn}, 62.3, all we have to show is that the diagonal functor 
$diag$: $\Delta \ra \Delta \times \Delta$
is initial.
For $(\underline{k_1},\underline{k_2}) \in Ob(\Delta \times \Delta)$, we have to show that the nerve of the over category $(diag \downarrow (\underline{k_1},\underline{k_2}))$ is contractible.

If $K$ is a simplicial set, $\Delta K$  (the ``barycentric subdivision'') denotes the category with $Ob \Delta K$ = simplicial set maps of the form $\Delta[n] \ra K$, and $Hom_{\Delta K}(a, b)$ = commutative diagrams of the form
\begin{center}
$\xymatrix@=.2in{
	\Delta[n] \ar[r]^a \ar[dr] &
	K \\
	& \Delta[m] \ar[u]_b
}$
\end{center}

It is known that $\Delta K$ has the same homotopy type as $K$.

We have an isomorphism of categories
\begin{center}
$(diag \downarrow (\underline{k_1},\underline{k_2})) = \Delta(\Delta[k_1] \times \Delta[k_2])$
\end{center}

It follows that $(diag \downarrow (\underline{k_1},\underline{k_2}))$ has the homotopy type of $\Delta[k_1] \times \Delta[k_2]$, therefore it is contractible.
\end{proof}

\bigskip
\bigskip

\bigskip
\bigskip
\bigskip

                % Local Variables: 
                % mode: latex 
                % tex-command: amslatex-command
                % tex-dvi-view-command: tex-oneside-view-command
                % tex-start-of-header: "\\special"
                % tex-end-of-header: "\\maketitle"
                % abbrev-mode: t
                % End:

%% file: appc.tex
\chapter{Edgewise Subdivision of Simplicial Diagrams}
\label{app3}

\bigskip
\bigskip

The material in this Appendix is a rewrite for model categories of the edgewise subdivision techniques of B\"ockstedt, Hsiang and Madsen \cite{BockstedtHsiangMadsen}.

The edgewise subdivision functor $F_k$ is defined for $k \geq 1$ by
\begin{center}
$ F_k : \Delta \lra \Delta $

$ F_k(\underline{n-1}) = \underline{kn - 1} $

\bigskip

\end{center} 
where if $\phi$ : $\underline{n_1-1} \lra \underline{n_2-1}$, then $F_k(\phi)$ : $\underline{kn_1 - 1} \lra \underline{kn_2 - 1}$ is given by 
\begin{center}
\[ F_k(\phi)(tn_1+i) = tn_2 + \phi(i) \phantom{fooo} (0 \leq t \leq k-1, 0 \leq i \leq n_1-1) \]
\end{center}

\bigskip

It is not hard to observe that $F_{kk'} = F_k F_{k'}$ and $F_1 = id_\Delta$.

$F_k$ is an initial functor, for $k \geq 1$ (see Prop. \ref{propFkInit}).

\bigskip

For $1 \leq l \leq k$, one constructs natural maps $u_k^l$ : $id_\Delta \lra F_k$ by 
\begin{center}
\bigskip
$u_k^l$ : $\underline{n-1} \lra F_k(\underline{n-1}) \phantom{fooo}$ given by

\bigskip

$u_k^l(i) = (l-1)n + i \phantom{fooo} (0 \leq i \leq n-1)$

\bigskip
\end{center} 

For fixed $k \geq 1$, the natural maps $u_k^l$ : $id_\Delta \lra F_k$ ($1 \leq l \leq k$) are pairwise {\it homotopic} maps (see Prop. \ref{proppsikl}), in the sense that, for any complete category ${\cat C}$ and for any cosimplicial object $X \tdot \in Ob{\cat C}^\Delta$, the $k$ cosimplicial maps
\begin{center}
\bigskip

$u_k^{l*}$ : $X \tdot \lra F_k^* X \tdot \phantom{fooo} (1 \leq l \leq k)$

\bigskip
\end{center} 
are pairwise homotopic: $u_k^{l*} \simeq u_k^{l'*}$ for $1 \leq l$, $l' \leq k$.

\bigskip

Our practical reason for writing this Appendix is Prop. \ref{proppsikl}, needed for proving that the cocycle relations in the proof of Chap. \ref{chap1}, Thm. \ref{thmBCdef} are satisfied.

\bigskip
\bigskip
\bigskip

\section{Edgewise Subdivision}

\bigskip
\bigskip

\subsection{Properties of ${\bf R}\lim$}
In this Section we will base our proofs on the following two results about homotopy limits:
\begin{lemma}[Diagram naturality of ${\bf R}\lim$]
\label{lemmaNatRlim}
Let ${\cat M}$ be a model category.
Let ${\cat D}$ and ${\cat D}'$ be small categories, and suppose we have two functors $U, V$ : ${\cat D} \xymatrix{\ar@<.5ex>[r] \ar@<-.5ex>[r] &} {\cat D}'$ and a natural map $u$ : $U \lra V$.

For any ${\cat D}'$-diagram in ${\cat M}$ denoted $X \tdot \in Ob({\cat M}^{{\cat D}'})$ the following diagram in ${\bf ho}{\cat M}$ is commutative, where the vertical maps are the obvious ones:
\begin{center}
\bigskip
$\xymatrix{
	{\bf R}\lim^{{\cat D}'} U^* X \tdot \ar[r]^{{\bf R}\lim^u} &
	{\bf R}\lim^{{\cat D}'} V^* X \tdot \\
	& {\bf R}\lim^{\cat D} X \tdot \ar[lu] \ar[u]
}$
\bigskip
\end{center}

Furthermore, given a natural map $v$ : $V \lra W$ then ${\bf R}\lim^v \comp {\bf R}\lim^u$ = ${\bf R}\lim^{vu}$.
If $id$ : $U \lra U$ is the identity natural map, then ${\bf R}\lim^{id}$ = $id$.
\end{lemma}

\begin{proof}
Easy consequence of the right adjointness property of ${\bf R}\lim$, described in Chap. \ref{chapHSS}, Sec. \ref{subsecHomotLim}.
\end{proof}

There is a dual for ${\bf L}\colim$ of Lemma \ref{lemmaNatRlim}, that we will not state for brevity.

\begin{thm}[Hirschhorn, \cite{Hirschhorn}, Thm. 20.6.11]
\label{thmHirsch}
Let ${\cat D}$ and ${\cat D}'$ be small categories, and suppose we have a functor $U$ : ${\cat D} \lra {\cat D}'$.
Then $U$ is initial if and only if for every model category ${\cat M}$ the natural map in ${\bf ho}{\cat M}$
\begin{center}
\bigskip
${\bf R}\lim^{\cat D} X \tdot \lra {\bf R}\lim^{{\cat D}'} U^* X \tdot$
\bigskip
\end{center}
is an isomorphism. 
\end{thm}
Again, we will not formulate the dual of Theorem \ref{thmHirsch}.

\bigskip
\bigskip

\subsection{Unfolded analogs of $F_k$ and $u_k^l$}

Consider the functor $G_k$ : $\Delta^{\times k} \lra \Delta$ from the $k$-fold product of $\Delta$ with itself, where $G$ is defined on objects by
\begin{center}
\bigskip
$G_k(\underline{n_1-1}, \underline{n_2-1}, ... , \underline{n_k-1})$ = $\underline{n_1 + n_2 + ... + n_k - 1}$
\bigskip
\end{center}
and $G_k$ is defined on maps by $\phi_l$ : $\underline{n_l-1} \lra \underline{n_l^\prime-1}$ ($1 \leq l \leq k$) going to $G_k(\phi_1, \phi_2, ... , \phi_k)$ : $\underline{n_1 + n_2 + ... + n_k - 1} \lra \underline{n_1^\prime + n_2^\prime + ... + n_k^\prime - 1}$ given by 
\begin{center}
\bigskip
$G_k(\phi_1, \phi_2, ... , \phi_k)(n_1 + n_2 + ... + n_{l-1} + i)$ = $n_1^\prime + n_2^\prime + ... + n_{l-1}^\prime + \phi_l(i) \phantom{foo} (1 \leq l \leq k, \phantom{f} 0 \leq i \leq n_l -1)$
\bigskip
\end{center}

\bigskip

Observe that $G_k$ is an ``unfolded'' version of $F_k$, in the sense that
\begin{center}
\bigskip
$F_k$ = $G_k \comp diag_k$ : $\xymatrix{\Delta \ar[r]^-{diag_k} & \Delta^{\times k} \ar[r]^-{G_k} & \Delta}$
\bigskip
\end{center}
where we denoted $diag_k$ : $\Delta \lra \Delta^{\times k}$ the $k$-fold diagonal.

\bigskip
\bigskip
\bigskip

We are next trying to create the right context to apply Lemma \ref{lemmaNatRlim} to the functor $G_k$ : $\Delta^{\times k} \lra \Delta$.
Fix any integer $l$, $1 \leq l \leq k$ and consider the projection on the $l$-th factor functor $\pi_k^l$ : $\Delta^{\times k} \lra \Delta$.

We construct a natural map $v_k^l$ : $\pi_k^l \lra G_k$, which should be thought of as an unfolded version of the natural map $u_k^l$ : $id_\Delta \lra F_k$.
The natural map $v_k^l$ is defined by
\begin{center}
\bigskip
$v_k^l$ : $\pi_k^l(\underline{n_1-1}, ... ,\underline{n_k-1})  \lra G_k(\underline{n_1-1}, ... ,\underline{n_k-1})$ 

\bigskip

$v_k^l$ : $\underline{n_l-1}  \lra \underline{n_1 + ... + n_k -1}$ 

\bigskip

$v_k^l(i) = n_1 + ... + n_{l-1} + i \phantom{fooo} (0 \leq i \leq n_l-1)$

\bigskip
\end{center} 

Observe that $v_k^l$ : $\pi_k^l \lra G_k$ applied to $diag_k$ is 
\begin{center}
$u_k^l$ : $\pi_k^l \comp diag_k = id_\Delta \lra G_k \comp diag_k = F_k$
\end{center}

\bigskip
\bigskip
\bigskip

We are ready to prove the essential Proposition that will lead to proving that the functor $F_k$ is initial.
The reader should compare the proof of this Proposition with the proof of Chap. \ref{chap1}, Thm. \ref{thmBCdef}.

\begin{prop}
\label{propGkinit0}
For any model category ${\cat M}$ and any cosimplicial object in ${\cat M}$ denoted $X \tdot \in Ob({\cat M}^\Delta)$, the natural map in ${\bf ho}{\cat M}$
\begin{center}
\[
{\bf R} {\lim}^{\Delta^{\times k}} \pi_k^{l*} X \tdot 
\xymatrix{\ar[r]^{{\bf R} {\lim}^{v_k^l}}_\cong &} 
{\bf R} {\lim}^{\Delta^{\times k}} G_k^* X \tdot
\]
\end{center}
is an isomorphism.
\end{prop}

\begin{proof}
Note that we can identify ${\bf R} {\lim}^\Delta X \tdot \xymatrix{\ar[r]^\cong &} {\bf R}{\lim}^{\Delta^{\times k}} \pi_k^{l*} X \tdot$ .

$G_k^* X \tdot$ is a $k$-fold cosimplicial object in ${\cat M}$.
The key observation is that we can succesively contract the first, second, ... (skip $l$-th), ... $k$-th cosimplicial dimensions of $G_k^* X \tdot$ .
The resulting object after contractions is just $X \tdot$ , and using Appendix \ref{app1}, Thm. \ref{thmHomot2} we get the desired isomorphism
\begin{center}
\[
{\bf R} {\lim}^\Delta X \tdot 
\xymatrix{\ar[r]^{{\bf R} {\lim}^{v_k^l}}_\cong &} 
{\bf R} {\lim}^{\Delta^{\times k}} G_k^* X \tdot
\]
\end{center}
\end{proof}

At this point, it is easy to prove:
\begin{prop}
\label{propGkinit}
The functor $G_k$ : $\Delta^{\times k} \lra \Delta$ is initial.
\end{prop}
\begin{proof}
Let ${\cat M}$ be any model category.
Choose {\bf any} $l$, $1 \leq l \leq k$. 

We use Lemma \ref{lemmaNatRlim} for ${\cat D}$ = $\Delta$, ${\cat D}'$ = $\Delta^{\times k}$, $U$ = $\pi_k^l$, $V$ = $G_k$, $u$ = $v_k^l$.
From the commutative diagram
\begin{center}
\bigskip
$\xymatrix{
	{\bf R}\lim^{\Delta^{\times k}} \pi_k^{l*} X \tdot \ar[r]^{{\bf R}\lim^{v_k^l}}_\cong &
	{\bf R}\lim^{\Delta^{\times k}} G_k^* X \tdot \\
	& {\bf R}\lim^\Delta X \tdot \ar[lu]^\cong \ar[u]
}$
\bigskip
\end{center}
we deduce that the natural map ${\bf R}\lim^\Delta X \tdot \lra {\bf R}\lim^{\Delta^{\times k}} G_k^* X \tdot$ is an isomorphism.

Since the model category ${\cat M}$ is arbitrary, from Thm. \ref{thmHirsch} we deduce that the functor $G_k$ is initial.
\end{proof}

\bigskip
\bigskip

\subsection{Proofs of properties of $F_k$ and $u_k^l$}

\begin{prop}
\label{propFkInit}
$F_k$ is an initial functor, for $k \geq 1$.
\end{prop}

\begin{proof}
Recall that $F_k$ = $G_k \comp diag_k$.

The functor $G_k$ is initial by Prop. \ref{propGkinit}.
The $k$-fold diagonal functor $diag_k$ is initial, using arguments developed in Appendix \ref{app2}.
By the transitivity property of initial functors, it follows that $F_k$ is initial.
\end{proof}

\bigskip
\bigskip

As a consequence of Lemma \ref{lemmaNatRlim} and Prop. \ref{propFkInit}, it is not hard to see that given a model category ${\cat M}$ and $X \tdot \in Ob({\cat M}^\Delta)$ then all the $k$ maps below coincide:
\begin{center}
\[
{\bf R} {\lim}^\Delta X \tdot 
\xymatrix{\ar[r]^{{\bf R} {\lim}^{u_k^l}}_\cong &} 
{\bf R} {\lim}^\Delta F_k^* X \tdot
\phantom{foo} (1 \leq l \leq k)
\]
\end{center}

A refinement of this result is given by:

\begin{prop}
\label{proppsikl}
Let ${\cat C}$ be any complete category, and $X \tdot$ a cosimplicial object in ${\cat C}$.
Then the $k$ cosimplicial maps
\begin{center}
\bigskip

$u_k^{l*}$ : $X \tdot \lra F_k^* X \tdot \phantom{fooo} (1 \leq l \leq k)$

\bigskip
\end{center} 
are pairwise homotopic: $u_k^{l*} \simeq u_k^{l'*}$ for $1 \leq l$, $l' \leq k$.
\end{prop}
 
\begin{proof}
There exists a natural map
\begin{center}
\bigskip

$h_k$ : $X \tdot \lra (F_k^* X \tdot)^{\Delta[k-1]}$

\bigskip

$h_k$ : $X^{n-1} \lra \times_{\Delta[k-1]_{n-1}} X^{kn-1} \phantom{fooo} (1 \leq n)$

\bigskip
\end{center} 
given on the factor in the product $\times_{\Delta[k-1]_{n-1}}$ corresponding to $\Phi$ : $\underline{n-1} \lra \underline{k-1}$ by the cosimplicial structure map $\Psi$ : $X^{n-1} \lra X^{kn-1}$ described as
\begin{center}
$\Psi$ : $\underline{n-1} \lra \underline{kn-1}$

\bigskip

$\Psi(i)$ = $\Phi(i)n + i \phantom{fooo} (0 \leq i \leq n-1)$ 
\end{center} 

Using the natural map $h_k$, it is easy to construct homotopies $u_k^{l*} \simeq u_k^{l'*}$ for $1 \leq l$, $l' \leq k$.
\end{proof}

\bigskip
\bigskip
\bigskip

                % Local Variables: 
                % mode: latex 
                % tex-command: amslatex-command
                % tex-dvi-view-command: tex-oneside-view-command
                % tex-start-of-header: "\\special"
                % tex-end-of-header: "\\maketitle"
                % abbrev-mode: t
                % End:

%% file: main.bbl
\begin{thebibliography}{99}

\bibitem{Andre1} M. Andr\'{e}: \emph{Homologie des alg\`{e}bres commutatives}. In: Die Grundlehren der mathematishen Wissenshaften in Einzerdarstellungen Band 206, Springer-Verlag, 1974.
%\bibitem{Avramov1} L. Avramov: \emph{Locally complete intersection homomorphisms and a conjecture of Quillen on the vanishing of cotangent cohomology}. Preprint 1997, 
%
%at http://www.math.purdue.edu/$\sim$avramov/papers.html
%\bibitem{Avramov2} L. Avramov: \emph{Infinite free resolutions}. Preprint 1996,
%
%at http://www.math.purdue.edu/$\sim$avramov/papers.html
%\bibitem{BDH} G. Baumslag, E. Dyer and A. Heller: \emph{The topology of discrete groups}. J. Pure Appl. Alg. 16 (1980), 1-47.
\bibitem{BockstedtHsiangMadsen}M. B\"ockstedt, W. C. Hsiang, I. Madsen: \emph{The cyclotomic trace and algebraic K-theory}. Invent. math. 111, 465-540 (1993).
\bibitem{Bousfield-Kan1} A. K. Bousfield, D. M. Kan: \emph{Homotopy limits, completions and localisations}. Lecture Notes in Mathematics 304, Springer, 1972.
\bibitem{Dwyer-Kan-Hirschhorn} W. G. Dwyer, D. M. Kan, P. S. Hirschhorn: \emph{Model categories and more general abstract homotopy theory: a work in what we like to think of as progress}. Preprint, draft March 28, 1997.

at http://www-math.mit.edu/$\sim$psh/kanmain.dvi.gz
\bibitem{Dwyer-Kan1} W. G. Dwyer, D. M. Kan: \emph{Simplicial localizations of categories}. J. Pure Appl. Alg. 17 (1980), 267-284.
\bibitem{Dwyer-Kan2} W. G. Dwyer, D. M. Kan: \emph{Calculating simplicial localizations}. J. Pure Appl. Alg. 18 (1980), 17-35.
\bibitem{Dwyer-Kan3} W. G. Dwyer, D. M. Kan: \emph{Function complexes in homotopical algebra}. Topology 19 (1980), 427-440.
%\bibitem{Goerss1} P. F. Goerss: \emph{On the Andr\'{e}-Quillen cohomology of commutative ${\mathbb F}_2$- algebras}. Ast\'{e}risque, vol. 186, Soc. Math. France, Paris 1990.
%\bibitem{Goerss2} P. F. Goerss: \emph{A Hilton-Milnor theorem for categories of simplicial algebras}. Am. J. Math. 111 (1989), 927-971.
\bibitem{Goerss-Jardine} P. J. Goerss, J. F. Jardine: \emph{Simplicial Homotopy Theory}. Preprint, version July 30, 1997.  

at http://www.math.uwo.ca/$\sim$jardine/papers/simp-sets/
\bibitem{Hirschhorn} P. S. Hirschhorn: \emph{Localization of model categories}. Preprint, draft April 27, 1998.

at http://www-math.mit.edu/$\sim$psh/
\bibitem{Hovey} M. Hovey: \emph{Model Categories}. Mathematical Surveys and Monographs No. 63, A.M.S. 1998.

%\bibitem{Kan-Thurston} D. Kan, W. Thurston: \emph{Every connected space has the homology of a} $K(\pi, 1)$. Topology 15 (1976), 253-258. 
\bibitem{MacLane} S. Mac Lane: \emph{Categories for the Working Mathematician}. Springer-Verlag New York, Heilderberg, Berlin 1971.
%\bibitem{Maunder} C. R. F. Maunder: \emph{A short proof of a theorem of Kan and Thurston}. Bull. London Math. Soc. 13 (1981), 325-327.
\bibitem{Meyer} J.-P. Meyer: \emph{Cosimplicial Homotopies}. Proc. A.M.S., Vol. 108, No. 1, pag 9-17.
%\bibitem{Miller1} H. Miller: \emph{The Sullivan conjecture on maps from classifying spaces}. Ann. Math, 120 (1984), 39-87. Correction, Ann. Math, 121 (1985), 605-609.
%\bibitem{Turner} J. Turner: \emph{On simplicial commutative rings with vanishing Andr\'{e}-Quillen homology}. Preprint, 1997, on the Hopf topology preprint archive.
\bibitem{Quillen1} D. Quillen: \emph{Homotopical algebra}. Lecture Notes in Mathematics 43. Heidelberg: Springer 1967.
\bibitem{Quillen2} D. Quillen: \emph{On the (co)homology of commutative rings}. Proc. Symp. Pure Math. 17 (1970), 65-87.
\bibitem{Quillen3} D. Quillen: \emph{On the homology of commutative rings}. Manuscript, M.I.T.





\end{thebibliography}
